%
%
%

%


\magnification=1200
\pretolerance=500 \tolerance=1000 \brokenpenalty=5000
\hsize=12.5cm   
\vsize=19cm
\hoffset=0.4cm
\voffset=1cm   
\parskip3pt plus 1pt
\parindent=0.6cm
\let\sl=\it


\font\seventeenbf=cmbx10 at 17.28pt

\font\twelvebf=cmbx10 at 12pt
\font\eightbf=cmbx8
\font\sixbf=cmbx6

\font\eighti=cmmi8
\font\sixi=cmmi6

\font\twelverm=cmr10 at 12pt
\font\eightrm=cmr8
\font\sixrm=cmr6

\font\eightsy=cmsy8
\font\sixsy=cmsy6

\font\eightit=cmti8
\font\eighttt=cmtt8
\font\eightsl=cmsl8

\font\seventeenbsy=cmbsy10 at 17.28pt

\font\twelvebsy=cmbsy10 at 12pt
\font\tenbsy=cmbsy10
\font\eightbsy=cmbsy8
\font\sevenbsy=cmbsy7
\font\sixbsy=cmbsy6
\font\fivebsy=cmbsy5

\font\tenmsa=msam10

\font\sevenmsa=msam7
\font\fivemsa=msam5
\newfam\msafam
  \textfont\msafam=\tenmsa
  \scriptfont\msafam=\sevenmsa
  \scriptscriptfont\msafam=\fivemsa

\font\tenmsb=msbm10
\font\eightmsb=msbm8
\font\sevenmsb=msbm7
\font\fivemsb=msbm5
\newfam\msbfam
  \textfont\msbfam=\tenmsb
  \scriptfont\msbfam=\sevenmsb
  \scriptscriptfont\msbfam=\fivemsb
\def\Bbb{\fam\msbfam\tenmsb}

\font\tenCal=eusm10
\font\sevenCal=eusm7
\font\fiveCal=eusm5
\newfam\Calfam
  \textfont\Calfam=\tenCal
  \scriptfont\Calfam=\sevenCal
  \scriptscriptfont\Calfam=\fiveCal
\def\Cal{\fam\Calfam\tenCal}

\font\teneuf=eusm10
\font\teneuf=eufm10
\font\seveneuf=eufm7
\font\fiveeuf=eufm5
\newfam\euffam
  \textfont\euffam=\teneuf
  \scriptfont\euffam=\seveneuf
  \scriptscriptfont\euffam=\fiveeuf
\def\euf{\fam\euffam\teneuf}

\font\seventeenbfit=cmmib10 at 17.28pt

\font\twelvebfit=cmmib10 at 12pt
\font\tenbfit=cmmib10
\font\eightbfit=cmmib8
\font\sevenbfit=cmmib7
\font\sixbfit=cmmib6
\font\fivebfit=cmmib5
\newfam\bfitfam
  \textfont\bfitfam=\tenbfit
  \scriptfont\bfitfam=\sevenbfit
  \scriptscriptfont\bfitfam=\fivebfit


\catcode`\@=11
\def\eightpoint{%
  \textfont0=\eightrm \scriptfont0=\sixrm \scriptscriptfont0=\fiverm
  \def\rm{\fam\z@\eightrm}%
  \textfont1=\eighti \scriptfont1=\sixi \scriptscriptfont1=\fivei
  \def\oldstyle{\fam\@ne\eighti}%
  \textfont2=\eightsy \scriptfont2=\sixsy \scriptscriptfont2=\fivesy
  \textfont\itfam=\eightit
  \def\it{\fam\itfam\eightit}%
  \textfont\slfam=\eightsl
  \def\sl{\fam\slfam\eightsl}%
  \textfont\bffam=\eightbf \scriptfont\bffam=\sixbf
  \scriptscriptfont\bffam=\fivebf
  \def\bf{\fam\bffam\eightbf}%
  \textfont\ttfam=\eighttt
  \def\tt{\fam\ttfam\eighttt}%
  \textfont\msbfam=\eightmsb
  \def\Bbb{\fam\msbfam\eightmsb}%
  \abovedisplayskip=9pt plus 2pt minus 6pt
  \abovedisplayshortskip=0pt plus 2pt
  \belowdisplayskip=9pt plus 2pt minus 6pt
  \belowdisplayshortskip=5pt plus 2pt minus 3pt
  \smallskipamount=2pt plus 1pt minus 1pt
  \medskipamount=4pt plus 2pt minus 1pt
  \bigskipamount=9pt plus 3pt minus 3pt
  \normalbaselineskip=9pt
  \setbox\strutbox=\hbox{\vrule height7pt depth2pt width0pt}%
  \let\bigf@ntpc=\eightrm \let\smallf@ntpc=\sixrm
  \normalbaselines\rm}
\catcode`\@=12

\def\eightpointbf{%
 \textfont0=\eightbf   \scriptfont0=\sixbf   \scriptscriptfont0=\fivebf
 \textfont1=\eightbfit \scriptfont1=\sixbfit \scriptscriptfont1=\fivebfit
 \textfont2=\eightbsy  \scriptfont2=\sixbsy  \scriptscriptfont2=\fivebsy
 \eightbf
 \baselineskip=10pt}

\def\tenpointbf{%
 \textfont0=\tenbf   \scriptfont0=\sevenbf   \scriptscriptfont0=\fivebf
 \textfont1=\tenbfit \scriptfont1=\sevenbfit \scriptscriptfont1=\fivebfit
 \textfont2=\tenbsy  \scriptfont2=\sevenbsy  \scriptscriptfont2=\fivebsy
 \tenbf}
        
\def\twelvepointbf{%
 \textfont0=\twelvebf   \scriptfont0=\eightbf   \scriptscriptfont0=\sixbf
 \textfont1=\twelvebfit \scriptfont1=\eightbfit \scriptscriptfont1=\sixbfit
 \textfont2=\twelvebsy  \scriptfont2=\eightbsy  \scriptscriptfont2=\sixbsy
 \twelvebf
 \baselineskip=14.4pt}

\def\seventeenpointbf{%
 \textfont0=\seventeenbf  \scriptfont0=\twelvebf  \scriptscriptfont0=\eightbf
 \textfont1=\seventeenbfit\scriptfont1=\twelvebfit\scriptscriptfont1=\eightbfit
 \textfont2=\seventeenbsy \scriptfont2=\twelvebsy \scriptscriptfont2=\eightbsy
 \seventeenbf
 \baselineskip=20.736pt}
 

\newdimen\srdim \srdim=\hsize
\newdimen\irdim \irdim=\hsize
\def\NOSECTREF#1{\noindent\hbox to \srdim{\null\dotfill ???(#1)}}
\def\SECTREF#1{\noindent\hbox to \srdim{\csname REF\romannumeral#1\endcsname}}
\def\INDREF#1{\noindent\hbox to \irdim{\csname IND\romannumeral#1\endcsname}}
\let\sectref=\SECTREF 
        
\newbox\titlebox   \setbox\titlebox\hbox{\hfil}
\newbox\sectionbox \setbox\sectionbox\hbox{\hfil}
\def\folio{\ifnum\pageno=1 \hfil \else \ifodd\pageno
           \hfil {\eightpoint\copy\sectionbox\kern8mm\number\pageno}\else
           {\eightpoint\number\pageno\kern8mm\copy\titlebox}\hfil \fi\fi}
\footline={\hfil}
\headline={\folio}

\def\titlerunning#1{\setbox\titlebox\hbox{\eightpoint #1}}
\def\title#1{\noindent\hfil$\smash{\hbox{\seventeenpointbf #1}}$\hfil
             \titlerunning{#1}\medskip}

\newcount\numbersection \numbersection=-1
\def\sectionrunning#1{\setbox\sectionbox\hbox{\eightpoint #1}
  \immediate\write1{\string\def \string\REF 
      \romannumeral\numbersection \string{%
      \noexpand#1 \string\dotfill \space \number\pageno \string}}}
\def\section#1{%
  \par\vskip0.5cm\penalty -100
  \vbox{\baselineskip=14.4pt\noindent{{\twelvepointbf #1}}}
  \vskip2pt
  \penalty 500
  \advance\numbersection by 1
  \sectionrunning{#1}}

\def\subsection#1|{%
  \par\vskip0.5cm\penalty -100
  \vbox{\noindent{{\tenpointbf #1}}}
  \vskip1pt
  \penalty 500}

\newcount\numberindex \numberindex=0  
\def\index#1#2{%
  \advance\numberindex by 1
  \immediate\write1{\string\def \string\IND #1%
     \romannumeral\numberindex \string{%
     \noexpand#2 \string\dotfill \space \string\S \number\numbersection, 
     p.\string\ \space\number\pageno \string}}}

\newdimen\itemindent \itemindent=\parindent

\def\item#1{\par\noindent\hangindent\itemindent%
            \rlap{#1}\kern\itemindent\ignorespaces}
\def\itemitem#1{\par\noindent\hangindent2\itemindent%
            \kern\itemindent\rlap{#1}\kern\itemindent\ignorespaces}
\def\itemitemitem#1{\par\noindent\hangindent3\itemindent%
            \kern2\itemindent\rlap{#1}\kern\itemindent\ignorespaces}

\def\qtq#1{\quad\text{#1}\quad}

\long\def\claim#1|#2\endclaim{\par\vskip 5pt\noindent 
{\tenpointbf #1.}\ {\sl #2}\par\vskip 5pt}

\def\proof{\noindent{\sl Proof}}

\def\today{\ifcase\month\or
January\or February\or March\or April\or May\or June\or July\or August\or
September\or October\or November\or December\fi \space\number\day,
\number\year}

\catcode`\@=11
\newcount\@tempcnta \newcount\@tempcntb 
\def\timeofday{{%
\@tempcnta=\time \divide\@tempcnta by 60 \@tempcntb=\@tempcnta
\multiply\@tempcntb by -60 \advance\@tempcntb by \time
\ifnum\@tempcntb > 9 \number\@tempcnta:\number\@tempcntb
  \else\number\@tempcnta:0\number\@tempcntb\fi}}
\catcode`\@=12

\def\bibitem#1&#2&#3&#4&%
{\hangindent=1.66cm\hangafter=1
\noindent\rlap{\hbox{\eightpointbf #1}}\kern1.66cm{\rm #2}{\sl #3}{\rm #4.}} 


\def\bC{{\Bbb C}}

\def\bN{{\Bbb N}}
\def\bP{{\Bbb P}}
\def\bQ{{\Bbb Q}}
\def\bR{{\Bbb R}}

\def\gm{{\euf m}}


\def\cC{{\Cal C}}
\def\cH{{\Cal H}}
\def\cI{{\Cal I}}
\def\cJ{{\Cal J}}
\def\cO{{\Cal O}}
\def\cP{{\Cal P}}
\def\cS{{\Cal S}}
\def\cT{{\Cal T}}
\def\cZ{{\Cal Z}}

\def\ampersand{\&}
\def\ld{,\ldots,}
\def\bu{{\scriptstyle\bullet}}

\def\square{{\hfill \hbox{
\vrule height 1.453ex  width 0.093ex  depth 0ex
\vrule height 1.5ex  width 1.3ex  depth -1.407ex\kern-0.1ex
\vrule height 1.453ex  width 0.093ex  depth 0ex\kern-1.35ex
\vrule height 0.093ex  width 1.3ex  depth 0ex}}}
\def\qed{\phantom{$\quad$}\hfill$\square$\medskip}
\def\hexnbr#1{\ifnum#1<10 \number#1\else
 \ifnum#1=10 A\else\ifnum#1=11 B\else\ifnum#1=12 C\else
 \ifnum#1=13 D\else\ifnum#1=14 E\else\ifnum#1=15 F\fi\fi\fi\fi\fi\fi\fi}
\def\msatype{\hexnbr\msafam}
\def\msbtype{\hexnbr\msbfam}
\mathchardef\restriction="3\msatype16   
\mathchardef\compact="3\msatype62
\mathchardef\smallsetminus="2\msbtype72   \let\ssm\smallsetminus
\mathchardef\subsetneq="3\msbtype28
\mathchardef\supsetneq="3\msbtype29
\mathchardef\leqslant="3\msatype36   \let\le\leqslant
\mathchardef\geqslant="3\msatype3E   \let\ge\geqslant
\mathchardef\ltimes="2\msbtype6E
\mathchardef\rtimes="2\msbtype6F

\let\ol=\overline

\let\wt=\widetilde

\let\text=\hbox
\def\build#1|#2|#3|{\mathrel{\mathop{\null#1}\limits^{#2}_{#3}}}


\def\Tr{\mathop{\rm Tr}\nolimits}

\def\Aut{\mathop{\rm Aut}\nolimits}
\def\GL{\mathop{\rm GL}\nolimits}

\def\Pic{\mathop{\rm Pic}\nolimits}

\def\Ricci{\mathop{\rm Ricci}\nolimits}

\def\codim{\mathop{\rm codim}\nolimits}

\def\div{\mathop{\rm div}\nolimits}

\def\ord{\mathop{\rm ord}\nolimits}
 
\def\diam{\mathop{\rm diam}\nolimits}
\def\pr{\mathop{\rm pr}\nolimits}
\def\dbar{{\overline\partial}}
\def\ddbar{{\partial\overline\partial}}


\def\reg{{\rm reg}}

\def\loc{{\rm loc}}

\def\ii{{\rm i}}
\def\Const{{\rm Const}}

\long\def\InsertFig#1 #2 #3 #4\EndFig{\par
\hbox{\hskip #1mm$\vbox to#2mm{\vfil\special{" 
(/home/demailly/psinputs/grlib.ps) run
#3}}#4$}}
\long\def\LabelTeX#1 #2 #3\ELTX{\rlap{\kern#1mm\raise#2mm\hbox{#3}}}

\def\RGBColor#1#2{#2}



\title{Semi-continuity of complex singularity}
\title{exponents and K\"ahler-Einstein}
\title{metrics on Fano orbifolds}
\medskip

\centerline{\twelvebf Jean-Pierre Demailly}\medskip
\centerline{\twelverm Universit\'e de Grenoble I, Institut Fourier}
\bigskip
\centerline{\twelvebf J\'anos Koll\'ar}\medskip
\centerline{\twelverm Princeton University, Department of Mathematics}

\titlerunning{J.-P.\ Demailly, J.\ Koll\'ar, 
  Semi-continuity of complex singularity exponents}

\vskip50pt
\noindent{\bf Abstract.} 
We introduce complex singularity exponents of plurisubharmonic functions
and prove a general semi-continuity result for them. This concept contains 
as a special case several similar concepts which have been considered e.g.
by Arnold and Varchenko, mostly for the study of hypersurface singularities. 
The plurisubharmonic version is somehow based on a reduction to the 
algebraic case, but it also takes into account more quantitative 
informations of great interest for complex analysis and complex
differential geometry.  We give as an application a new derivation of
criteria for the existence of K\"ahler-Einstein metrics on certain
Fano orbifolds, following Nadel's original ideas (but with a drastic
simplication in the technique, once the semi-continuity result is taken
for granted). In this way, 3 new examples of rigid K\"ahler-Einstein
Del Pezzo surfaces with quotient singularities are obtained.
\vskip20pt

\noindent{\bf R\'esum\'e.}
Nous introduisons les exposants de singularit\'es complexes des
fonctions plurisousharmoniques et d\'emontrons un th\'eor\`eme de
semi-continuit\'e g\'en\'eral pour ceux-ci. Le concept \'etudi\'e
contient comme cas particulier des concepts voisins qui ont \'et\'e
consid\'er\'es par exemple par Arnold et Varchenko, principalement
pour l'\'etude des singularit\'es d'hypersurfaces.  La version
plurisousharmonique repose en d\'efinitive sur une r\'eduction au cas
alg\'ebrique, mais elle prend aussi en compte des informations
quantitatives d'un grand int\'er\^et pour l'analyse complexe et la
g\'eom\'etrie diff\'erentielle complexe.  Nous d\'ecrivons en
application une nouvelle approche des crit\`eres d'existence de
m\'etriques K\"ahler-Einstein pour les vari\'et\'es de Fano, en nous
inspirant des id\'ees originales de Nadel -- mais avec des
simplifications importantes de la technique, une fois que le
r\'esultat de semi-continuit\'e est utilis\'e comme outil de
base. Gr\^ace \`a ces crit\`eres, nous obtenons trois nouveaux
exemples de surfaces de Del Pezzo \`a singularit\'es quotients,
rigides, poss\'edant une m\'etrique de K\"ahler-Einstein.
\vfill\break

\noindent{\bf Contents}
\vskip5pt

{\eightpoint
\def\REF{\S 0. Introduction \dotfill 2}

\sectref{0}
\sectref{1}
\sectref{2}
\sectref{3}
\sectref{4}
\sectref{5}
\sectref{6}}
\vskip8pt

\noindent
{\bf Key words:} complex singularity exponent, Arnold multiplicity,
semicontinuity property, log resolution, log canonical singularities,
effective divisor, coherent ideal sheaf, adjunction theory, inversion
of adjunction, plurisubharmonic function, multiplier ideal sheaf,
closed positive current, Lelong number, Monge-Amp\`ere equation,
Aubin-Calabi-Yau theorem, Fano variety, orbifold
\vskip10pt
\noindent
{\bf A.M.S. Classification 1991:} 14B05, 14J45, 32C17, 32J25, 32S05
\vskip20pt

\section{\S0. Introduction}

The purpose of this work is to show how complex analytic methods (and
more specifically $L^2$ estimates for $\dbar$) can provide effective
forms of results related to the study of complex singularities.
We prove in particular a strong form of the semi-continuity theorem for
``complex singularity exponents'' of plurisubharmonic (psh) functions.
An application to the existence of K\"ahler-Einstein metrics on certain
Fano orbifolds will finally be given as an illustration of this result.

We introduce the following definition as a quantitative way of
measuring singu\-la\-rities of a psh function $\varphi$ (the basic
definition even makes sense for an arbitrary measurable
function~$\varphi$, though it is unlikely to have any good properties
in that case). Our approach is to look at the $L^1$ integrability of
$\exp(-2c\varphi)$ in terms of the Lebesgue measure in some local
coordinates. Several other types of analytic or algebraic objects
(holomorphic functions, coherent ideal sheaves, divisors, currents,
etc) can be treated as special cases of this formalism.

\claim 0.1. Definition|Let $X$ be a complex manifold and $\varphi$ be
a plurisubharmonic $($psh$)$ function on~$X$. For any compact set
$K\subset X$, we introduce the ``complex singularity exponent'' of
$\varphi$ on $K$ to be the nonnegative number
\index{CO}{Complex singularity exponent}
$$
c_K(\varphi)=\sup\big\{c\ge 0\,:\,\exp(-2c\varphi)
\text{ is $L^1$ on a neighborhood of $K$}\big\},
$$
and we define the ``Arnold multiplicity'' to be
$\lambda_K(\varphi)=c_K(\varphi)^{-1}\,:$
$$
\lambda_K(\varphi)=\inf\big\{\lambda>0\,:\,\exp(-2\lambda^{-1}\varphi)
\text{ is $L^1$ on a neighborhood of $K$}\big\}.
$$
If $\varphi\equiv -\infty$ near some connected component of $K$, we put
of course $c_K(\varphi)=0$, $\lambda_K(\varphi)=+\infty$.
\endclaim

The singularity exponent $c_K(\varphi)$ only depends on the singularities
of $\varphi$, namely on the behavior of $\varphi$ near its $-\infty$
poles. Let $T$ be a closed positive current of bidegree $(1,1)$
on~$X$. Since $c_K(\varphi)$ remains unchanged if we replace $\varphi$
with $\psi$ such that $\psi-\varphi$ is bounded, we see that it is
legitimate to define
$$
c_K(T)=c_K(\varphi),\qquad\lambda_K(T)=\lambda_K(\varphi)
\leqno(0.1.1)
$$
whenever $\varphi$ is a (local) potential of $T$, i.e. a psh function
$\varphi$ such that $dd^c\varphi=T$, where $d^c=(2\pi\ii)^{-1}
(\partial-\dbar)$. In particular, if $D$ is an effective integral
divisor, we have $c_K([D])=c_K(\log|g|)$ where $[D]$ is the
current of integration over $D$ and $g$ is a (local) generator of
$\cO(-D)$. When $f$ is a holomorphic function, we write simply
$c_K(f)$, $\lambda_K(f)$ instead of $c_K(\log|f|)$, $\lambda_K(\log|f|)$.
For a coherent ideal sheaf $\cI=(g_1,\,\ldots\,,g_N)$
we define in a similar way $c_K=c_K(\log(|g_1|+\cdots+|g_N|))$. It is
well known that $c_K(f)$ is a rational number, equal to the largest
root of the Bernstein-Sato polynomial of $|f|^{2s}$ on a neighborhood
of~$K$ ([Lin89], see also [Kol97]); similarly $c_K(\cI)\in\bQ_+$ for
any coherent ideal sheaf. Our main result consists in the
following semi-continuity theorem.

\claim 0.2. Main Theorem|Let $X$ be a complex manifold. Let $\cZ^{1,1}_+(X)$
denote the space of closed positive currents of type $(1,1)$ on~$X$, equipped
with the weak topology, and let $\cP(X)$ be the set of locally $L^1$ psh 
functions on $X$, equipped with the topology of $L^1$ convergence on 
compact subsets $(=$~topology induced by the weak topology$)$. Then
\smallskip
\item{\rm(1)} The map $\varphi\mapsto c_K(\varphi)$ is lower semi-continuous 
on $\cP(X)$, and the map $T\mapsto c_K(T)$ is lower semi-continuous on
$\cZ^{1,1}_+(X)$.
\smallskip
\item{\rm(2)} {\rm(``Effective version'')}. Let $\varphi\in\cP(X)$ be given.
If $c<c_K(\varphi)$ and $\psi$ converges to $\varphi$ in $\cP(X)$, then
$e^{-2c\psi}$ converges to $e^{-2c\varphi}$ in $L^1$ norm over some
neighborhood $U$ of~$K$.
\smallskip\noindent
As a special case, one gets:
\item{\rm(3)} The map $\cO(X)\ni f\mapsto c_K(f)$ is
lower semi-continuous with respect to the topology of uniform convergence
on compact sets $($uniform convergence on a fixed neighborhood of~$K$ is
of course enough$)$. Moreover, if $c<c_K(f)$ and $g$ converges to~$f$ 
in~$\cO(X)$, then $|g|^{-2c}$ converges to $|f|^{-2c}$ in $L^1$ on some
neighborhood $U$ of~$K$.
\vskip0pt
\endclaim

In spite of their apparent simplicity, the above statements reflect
rather strong semi-continuity properties of complex singularities
under ``variation of parameters''. Such properties have been used
e.g.\ by Angehrn-Siu [AnSi95] in their approach of the Fujita conjecture,
and our arguments will borrow some of their techniques in section~\S$\,$3.

Theorem~0.2 is by nature a purely local result, which is easily seen
to be equivalent to the special case when $K=\{x\}$ is a single point
and $X$ is a small ball centered at~$x$. The proof is made in several
steps. The ``analytic part'' consists in a reduction of (1) and (2) to
(3), and in the proof of the effective estimates leading to the
convergence statements in (2) and (3) [by contrast, the qualitative
part of (3) can be obtained in a purely algebraic way]. The reduction
to the holomorphic case (3) is based on the fact that plurisubharmonic
functions can be very accurately approximated (both from the point of
view of singularities and of $L^1_\loc$ topology) by special functions
of the form
$$
\alpha\log(|g_1|+\cdots+|g_N|),\qquad \alpha\ge 0,\leqno(0.2.4)
$$
where the $g_j$ are holomorphic functions. The existence of
approximations as in (0.2.4) depends in an essential way on the
Ohsawa-Takegoshi $L^2$ extension theorem ([OhT87], [Ohs88]), see
[Dem92,~93] and \S2,~\S4. One is then reduced to the proof for a single
holomorphic function (that is, to a psh function of the form
$\log|f|$), by taking a suitable generic linear combination
$f=\sum\alpha_j g_j$. Another essential idea is to truncate the Taylor
expansion of $f$ at~$x$ at some order $k$. It can then be shown that
this affects $c_x(f)$ only by a perturbation that is under uniform
control. In fact, the singularity exponent $c_x(f)$ is subadditive
on holomorphic functions:
$$
c_x(f+g)\le c_x(f)+c_x(g),\qquad \forall f,g\in\cO_{X,x}.\leqno(0.2.5)
$$
If $p_k$ is the truncation at order $k$ of the Taylor series, one deduces
immediately from (0.2.5) that
$$
|c_x(f)-c_x(p_k)|\le {n\over k+1}.\leqno(0.2.6)
$$
In this way, the proof is reduced to the case of polynomials of
given degree. Such polynomials only depend on finitely many
coefficients, thus the remaining lower semi-continuity property to be
proved is that of the function $t\mapsto c_x(P_t)$ when $P_t$ is a
family of polynomials depending holomorphically on some
parameters~$t=(t_1\ld t_N)$.  This is indeed true, as was already
observed by Varchenko [Var82,~83].  An algebraic proof can be given by
using a log resolution of singularities with parameters. Here, however,
a special attention to effective estimates must be paid to prove the
convergence statements in (2) and (3). For instance, it is necessary to
get as well an effective version of (0.2.6); the Ohsawa-Takegoshi $L^2$
extension theorem is again crucial in that respect.

As a consequence of our main theorem, we give a more natural proof of
the results of Siu [Siu87,~88], Tian [Tia87] and Nadel [Nad89,~90] on 
the existence of K\"ahler-Einstein metric on Fano manifolds admitting 
a sufficiently big group of symmetries. The main
point is to have sufficient control on the ``multiplier ideal sheaves''
which do appear in case the K\"ahler-Einstein metric fails to exist. This
can be dealt with much more easily through our semi-continuity theorem,
along the lines suggested in Nadel's note [Nad89] (possibly because of the
lack of such semi-continuity results, the detailed version [Nad90] relies
instead on a rather complicated process based on a use of ``uniform''
$L^2$ estimates for sequences of Koszul complexes; all this disappears
here, thus providing a substantially shorter proof). We take the
opportunity to adapt Nadel's result to Fano orbifolds. This is mostly a
straightforward extension, except that we apply intersection
inequalities for currents rather than the existence of a big finite
group of automorphisms to derive sufficient criteria for the existence
of K\"ahler-Einstein metrics. In this way, we produce 3 new ``exotic
examples'' of rigid Del Pezzo surfaces with quotient singularities
which admit a K\"ahler-Einstein orbifold metric.

We would like to thank R.R.~Simha for useful discussions which got us
started with the idea of simplifying Nadel's approach. We also thank
Mongi Blel for sharing several viewpoints on the semicontinuity
properties of psh functions, and Jeff McNeal for pointing out a slight
inaccuracy in our original calculation of volumes of analytic tubes.

\section{\S1. Complex singularity exponent and Arnold multiplicity}

Let $X$ be a complex manifold and $\varphi$ a psh function of~$x$.
The concepts of ``complex singularity exponent'' $c_K(\varphi)$ and
``Arnold multiplicity'' $\lambda_K(\varphi)$ of $\varphi$ along a
compact set $K\subset X$ have been defined in~0.1. An equivalent
definition can be given in terms of asymptotic estimates for the
volume of sublevel sets $\{\varphi<\log r\}$.

\claim 1.1. Variant of the definition|Let $K\subset X$ be a compact set,
$U\compact X$ a relatively compact neighborhood of $K$, and let $\mu_U$ 
be the Riemannian measure on $U$ associated with some choice of 
hermitian metric $\omega$ on~$X$. Then
$$
c_K(\varphi)=\sup\big\{c\ge 0\,;\,r^{-2c}\mu_U(\{\varphi{<}\log r\})
\text{ is bounded as $r\to 0$, for some $U\supset K$}\big\}.
$$
\endclaim

The equivalence with the earlier Definition~0.1 follows immediately from
the elementary inequalities
$$
r^{-2c}\mu_U(\{\varphi<\log r\})\le\int_U e^{-2c\varphi}\,dV_\omega
\le\mu_U(U)+\int_0^1 2c\,r^{-2c}\mu_U(\{\varphi<\log r\})\,{dr\over r},
$$
A first important observation is that $c_K(\varphi)$ and $\lambda_K(\varphi)$
depend only on the local behavior of~$\varphi$:

\claim 1.2. Proposition|Given a point $x\in X$, we write $c_x(\varphi)$
instead of $c_{\{x\}}(\varphi)$. Then
$$
c_K(\varphi)=\inf_{x\in K}c_x(\varphi),\qquad
\lambda_K(\varphi)=\sup_{x\in K}\lambda_x(\varphi).
$$
\endclaim

The statement is clear from the Borel-Lebesgue Lemma. 
When $x$ is a pole, that is, when $\varphi(x)=-\infty$, the Arnold
multiplicity $\lambda_x(\varphi)$ actually measures the ``strength'' of
the singularity of $\varphi$ in a neighborhood of $x$. (It actually
``increases'' with the singularity, and if $x$ is not a pole, we have 
$c_x(\varphi)=+\infty$, $\lambda_x(\varphi)=0\,$; see Prop.~1.4 below.)
We now deal with various interesting special cases:

\claim 1.3. Notation|\smallskip
\item{\rm(1)} If $f$ is a holomorphic function on $X$, we set
$c_K(f)=c_K(\log|f|)$.
\smallskip
\item{\rm(2)} If $\cI\subset\cO_X$ is a coherent ideal sheaf, generated
by functions $(g_1\ld g_N)$ on a neighborhood of $K$, we put
$$c_K(\cI)=c_K\big(\log(|g_1|+\cdots+|g_N|)\big).$$
\smallskip
\item{\rm(3)} If $T$ is a closed positive current of bidegree $(1,1)$ on
$X$ which can be written as $T=dd^c\varphi$ on a neighborhood of $K$,
we set $c_K(T)=c_K(\varphi)$.
\smallskip
\noindent
$($If no global generators exist in {\rm(2)} or no global potential 
$\varphi$ exists in {\rm(3)}, we just split $K$ in finitely many pieces
and take the infimum, according to Prop.~$(1.2))$.
\smallskip
\item{\rm(4)} If $D$ is an effective divisor with rational or real 
coefficients, we set
$$
c_K(D)=c_K([D])=c_K(\cO(-D))=c_K(g)=c_K(\log|g|)
$$
where $D$ is the current of integration over $D$ and
$g$ is a local generator of the principal ideal sheaf $\cO(-D)$.
\vskip0pt
\endclaim

No confusion should arise from the above definitions, especially since 
$c_K(\cI)$ does not depend on the choice of generators of $\cI$. We use
similar conventions of notation for $\lambda_K(\varphi)$. The number
$$
c_x(f)=\sup\big\{c\,;\,|f|^{-2c}\text{ is $L^1$ on a neighborhood of $x$}
\big\}=\lambda_x(f)^{-1}
$$ 
is clearly a measure of the singularities of the hypersurface $\{f=0\}$ at 
point~$x$.  This number came up in the litterature many times under different
names. By [Lin89], $c_x(f)$ is the largest root of the Bernstein-Sato 
polynomial associated to the germ of $f$ around~$p$. If~$x$ is an isolated 
singularity of $\{f=0\}$, then $c_x(f)=\min\{1,\beta_{\bC}(f_x)\}$  where 
$\beta_{\bC}(f_x)$ is the complex singular index as defined in [ArGV84], 
vol.II, Sec.~13.1.5; the same thing is called ``complex singularity
exponent'' in [Var92]. See [Kol97] for a discussion of these
questions and for related results. 

\claim 1.4. Elementary properties|Let $\cI$, $\cJ$ be coherent ideals
on~$X$ and let $\varphi$, $\psi$ be psh functions. Denote by $x$ a point
in $X$ and let $K\subset X$ be a compact subset.
\smallskip
\item{\rm(1)} The function $x\mapsto c_x(\varphi)$ is lower
semi-continuous for the holomorphic Zariski topology$\,;$
\smallskip
\item{\rm(2)} If $\varphi\le\psi$, then $c_K(\varphi)\le c_K(\psi)\,;$
\hfill\break
If $\cI\subset\cJ$, then $c_K(\cI)\le c_K(\cJ)$.
\smallskip
\item{\rm(3)} 
$\lambda_K(\varphi+\psi)\le\lambda_K(\varphi)+\lambda_K(\psi)\,$;\hfill\break
$\lambda_K(\cI\cJ)\le\lambda_K(\cI)+\lambda_K(\cJ)$.
\smallskip
\item{\rm(4)} $\lambda_K(\alpha\varphi)=\alpha\,\lambda_K(\varphi)$
for all $\alpha\in\bR_+\,;$\hfill\break
$\lambda_K(\cI^m)=m\lambda_K(\cI)$ for all integers $m\in\bN$.
\smallskip
\item{\rm(5)} Let $\cI=(g_1\ld g_N)$ and let 
$$
\ol\cI=\big\{f\in\cO_{\Omega,x}\,,\,x\in\Omega\,;
\text{ $\exists C\ge 0$, $|f|\le C\max|g_j|$ near $x$}\big\}
$$ 
be the integral closure of~$\cI$. Then $c_K(\ol\cI)=c_K(\cI)$.
\smallskip
\item{\rm(6)} If the zero variety germ $V(\cI_x)$ contains a 
$p$-codimensional irreducible component, then $c_x(\cI)\le p$, i.e.
$\lambda_x(\cI)\ge 1/p$. 
\smallskip
\item{\rm(7)} If $\cI_Y$ is the ideal sheaf of a $p$-codimensional
subvariety $Y\subset\Omega$, then \hbox{$c_x(\cI_Y)=p$} at every
nonsingular point of~$Y$.
\item{\rm(8)} Define the vanishing order $\ord_x(\cI)$ of $\cI$ at
$x$ to be the supremum of all integers $k$ such that $\cI_x\subset\gm_x^k$,
where $\gm_x\subset\cO_x$ is the maximal ideal. Then
$$
{1\over n}\ord_x(\cI)\le\lambda_x(\cI)\le\ord_x(\cI).
$$
More generally, if $\nu_x(\varphi)$ is the Lelong number of $\varphi$ at~$x$,
then 
$$
{1\over n}\nu_x(\varphi)\le\lambda_x(\varphi)\le\nu_x(\varphi).
$$
\endclaim

\proof. (1) Fix a point $x_0$ and a relatively compact coordinate ball
$B:=B(x_0,r)\compact X$. For every $c\ge 0$, let $\cH_{c\varphi}(B)$ be 
the Hilbert space of holomorphic functions on $B$ with finite weighted 
$L^2$ norm
$$\Vert f\Vert_c^2=\int_B|f|^2e^{-2c\varphi}dV$$
where $dV$ is the Lebesgue volume element in $\bC^n$, $n=\dim_\bC X$.
A fundamental consequence of H\"ormander's $L^2$ estimates
(H\"ormander-Bombieri-Skoda theorem [H\"or66], [Bom70], [Sko75]) states
that there is an element $f\in\cH_{c\varphi}(B)$ with $f(x)=1$ whenever
$e^{-2c\varphi}$ is $L^1$ on a neighborhood of~$x$. Hence
$$
\big\{x\in B\,;\,c_x(\varphi)\le c_0\big\}\cap B=
\bigcap_{f\,\in\,\bigcup_{c>c_0}\cH_{c\varphi}(B)}f^{-1}(0)
$$
is an analytic set. This proves the holomorphic Zariski lower
semi-continuity.

All other properties are direct consequences of the definitions and do
not require ``hard'' analysis: (2), (4), (5) are immediate; (3) is a
consequence of the H\"older inequality; (6,7) follow from the fact that
the function $(\sum_{j\le p}|z_j|^2)^{-c}$ is locally integrable along
$z_1=\cdots=z_p=0$ if and only if $c<p\,$;
Finally, (8) is a well-known result of Skoda [Sko72], depending on the
basic properties of Lelong numbers and a use of standard kernel
techniques.\qed

In the case of an ideal sheaf, the following lemma reduces the computation 
of $c_x(\cI)$ to the case of a principal ideal (possibly after raising
$\cI$ to some power $\cI^m$).

\claim 1.5. Proposition|Let $(g_1\ld g_p)$ be holomorphic functions defined
on an open set $\Omega\subset\bC^n$ and let $x\in V(g_1\ld g_p)$. Then
$$c_x(\alpha_1g_1+\cdots+\alpha_pg_p)\le
\min\big\{c_x(g_1\ld g_p)\,,\,1\big\}$$
for all coefficients $(\alpha_1\ld\alpha_p)\in\bC^p$. Moreover, 
the equality occurs for all $(\alpha_1\ld\alpha_p)$ in the complement of a
set of measure zero in~$\bC^p$. In particular, if $\cI$ is an arbitrary
ideal and $c_x(\cI)\le 1$, there is a principal ideal $(f)\subset\cI$
such that $c_x(f)=c_x(\cI)$.
\endclaim

\proof. The inequality is obvious, since $c_x(\alpha_1g_1+\cdots+\alpha_pg_p)
\le 1$ by (1.4.6) on the one hand, and 
$$
\big|\alpha_1g_1+\cdots+\alpha_pg_p\big|^{-2c}
\ge\Big(\sum|\alpha_j|^2\Big)^{-c}\Big(\sum|g_j|^2\Big)^{-c}
$$
on the other hand. Now, fix $c<\min\{c_x(g_1\ld g_p)\,,\,1\}$.
There is a neighborhood $U_c$ of $x$ on which
$$
\leqalignno{
\int_{|\alpha|=1}d\sigma(\alpha)\int_{U_c}\big|\alpha_1g_1(z)
&{}+\cdots+\alpha_pg_p(z)\big|^{-2c}dV(z)\cr
&=A_c\int_{U_c}\Big(\sum|g_j(z)|^2\Big)^{-c}dV(z)<+\infty,&(1.5.1)\cr}
$$
where $d\sigma$ is the euclidean area measure on the unit sphere
$S^{2n-1}\subset\bC^n$ and $A_c>0$ is a constant. The above identity
follows from the formula
$$
\int_{|\alpha|=1}|\alpha\cdot w|^{-2c}d\sigma(\alpha)=A_c|w|^{-2c},
$$
which is obvious by homogeneity, and we have $A_c<+\infty$ for
$c<1$. The finiteness of the right hand side of (1.5.1) implies
that the left hand side is finite for all values $\alpha$ in the
complement $\bC^p\ssm N_c$ of a negligible set. Therefore
$c_x(\alpha_1g_1+\cdots+\alpha_pg_p)\ge c$, and by taking
the supremum over an increasing sequence of values $c_\nu$ converging
to~$\min\{c_x(g_1\ld g_p)\,,\,1\}$, we conclude that the equality holds 
in Proposition~1.5 for all $\alpha\in\bC^p\ssm\bigcup N_{c_\nu}$.\qed

\claim 1.6. Remark|{\rm It follows from Theorem~3.1 below that the
exceptional set of values $(\alpha_1\ld\alpha_p)$ occurring in
Prop.~1.5 is in fact a closed algebraic cone in $\bC^p$.\qed}
\endclaim

The singularity exponent $c_K(\cI)$ of a coherent ideal sheaf
$\cI\subset\cO_X$ can be computed by means of a ``log resolution''
of~$\cI$, that is, a composition $\mu:\wt X\to X$ of blow-ups with
smooth centers such that $\mu^\star\cI=\cO_{\wt X}(-D)$ is an
invertible sheaf associated with a normal crossing divisor $D$ in
$\wt X$ (such a log resolution always exists by Hironaka [Hir64]).
The following proposition is essentially well known (see e.g.\
[Kol95a] 10.7).

\claim 1.7. Proposition|Let $X$ be a complex manifold, $\cI\subset\cO_X$
a coherent ideal sheaf, and let $\mu:\wt X\to X$ be a modification
$(=$~proper bimeromorphic morphism$)$ such that
$\mu^\star\cI=\cO_{\wt X}(-D)$ is an invertible sheaf. Assume that $\wt X$
is normal and let $E_i\subset\wt X$ denote either an exceptional divisor
of $\mu$ or an irreducible component of~$D$. Write
$$
K_{\wt X}=\mu^\star K_X+\sum a_iE_i\qtq{and} D=\sum b_iE_i,
$$
where $a_i=0$ if $E_i$ is not a component of the exceptional divisor
of~$\mu$ $($resp.\ $b_i=0$ if $E_i$ is not a component of $D)$. Then:
\smallskip
\item{\rm(1)} 
$ 
\displaystyle c_K(\cI)\le \min_{i\,:\,\mu(E_i)\cap K
\neq \emptyset} \left\{{a_i+1\over b_i}\right\}.  
$
\smallskip
\item{\rm(2)} Equality holds if $\wt X$ is smooth and $\sum E_i$ is a
divisor with normal crossings.
\smallskip
\item{\rm(3)} If $g=(g_1,\,\ldots\,,g_N)$ are generators of $\cI$ in a
neighborhood of $K$, then for any sufficiently small neighborhood
$U$ of $K$ there is a volume estimate
$$
C_1r^{2c}\le\mu_U(\{|g|<r\})\le C_2 r^{2c}|\log r|^{n-1},
\qquad \forall r<r_0
$$
with $n=\dim_\bC X$, $c=c_K(\cI)$ and $C_1,\,C_2,\,r_0>0$.
\vskip0pt
\endclaim

\proof. Since the question is local, we may assume that $\cI$ is generated
by holomorphic functions $g_1,\,\ldots\,,g_N\in\cO(X)$. Then (1) and (2)
are straightforward consequences of the Jacobian formula for a 
change of variable: if $U$ is an open set in $X$, the change $z=\mu(\zeta)$
yields
$$
\int_{z\in U}|g(z)|^{-2c}dV(z)=
\int_{\zeta\in\mu^{-1}(U)}|g\circ\mu(\zeta)|^{-2c}|J_\mu(\zeta)|^2
d\wt V(\zeta)
$$
where $J_\mu$ is the Jacobian of $\mu$, and $dV$, $d\wt V$ are volume
elements of $X$, $\wt X$ respectively (embed $\wt X$ in some smooth
ambient space if necessary). Now, if $h_i$ is a generator of $\cO(-E_i)$
at a smooth point $\wt x\in \wt X$, the divisor of $J_\mu$ is by
definition $\sum a_iE_i$ and $\mu^\star\cI=\cO(-\sum b_iE_i)$.
Hence, up to multiplicative bounded factors,
$$
|J_\mu|^2\sim\prod|h_i|^{2a_i},\qquad
|g\circ\mu|^2\sim\prod|h_i|^{2b_i}\qquad\text{near $\wt x$,}
$$
and $|g\circ\mu|^{-2c}|J_\mu|^2$ is $L^1$ near $\wt x$ if and
only if $\prod|h_i|^{-2(cb_i-a_i)}$ is $L^1$. A necessary condition is
that $cb_i-a_i<1$ whenever $E_i\ni\wt x$. We therefore get the
necessary condition $c<\min_{i\,:\,\mu(E_i)\cap
K\ne\emptyset}\{(a_i+1)/b_i\}$, and this condition is necessary and
sufficient if $\sum E_i$ is a normal crossing divisor.

\noindent For (3), we choose $(\wt X,\cO(-D))$ to be a (nonsingular) log 
resolution of $\cI$. The volume $\mu_U(\{|g|<r\})$ is then given by 
integrals of the form
$$
\int_{\mu^{-1}(U)\cap\{\zeta\in\wt U_\alpha,\,\prod|h_i|^{b_i}<r\}}
\prod|h_i(\zeta)|^{2a_i}dV(\zeta)
\leqno(1.7.4)
$$
over suitable coordinate charts $\wt U_\alpha\subset\wt X$. An
appropriate change of variable $\zeta\mapsto w$, $w_i=h_i^{b_i}(\zeta)$,
$w_j=\zeta_{k_j}$ (where $i$ runs over the set of indices such that
$b_i>0$ and $j$ over a disjoint set of indices) and a use of a partition
of unity leads to estimate (1.7.4) by a linear combination of integrals
of the form
$$
\int_{P(r)}\prod|w_i|^{2(a_i+1)/b_i-2}dV(w)\qquad
\text{where~~$P(r)=\{\max|w_i|<1,\,\prod|w_i|<r\}$}
$$
(we assume here that a partial integration with respect to the $w_j$'s
has already been performed). The lower bound $C_1r^{2c}$ is obtained
by restricting the domain of integration to a neighborhood of a
point in the unit polydisk such that only one coordinate $w_i$
vanishes, precisely for $i$ equal to the index achieving the
minimum of $(a_i+1)/b_i$.
The upper bound $C_2r^{2c}|\log r|^{n-1}$, $c=\min(a_i+1)/b_i$,
is obtained by using the inequalities
$$
\eqalign{
\prod|&w_i|^{2(a_i+1)/b_i-2}\le\Big(\prod|w_i|\Big)^{2c-2}\le r^{2c-2}, 
\qquad\forall w\in P(r),\cr
\mu(P(r))&=\int_{\{\max(|w_1|,\ldots,|w_{n-1}|)<1\}}
\pi\,\min\Big({r^2\over |w_1|^2\cdots|w_{n-1}|^2},1\Big)
\prod_{i=1}^{n-1}dV(w_i)\cr
&\le\pi\int_{\{\exists i\,;\,|w_i|<r\}}
\prod_{i=1}^{n-1}dV(w_i)+\pi r^2
\int_{\{\forall i\,;\,r\le|w_i|<1\}}
\prod_{i=1}^{n-1}{dV(w_i)\over|w_i|^2}\cr
&\le C_2r^2|\log r|^{n-1}.\cr}
$$
It should be observed that
much finer estimates are known to exist; in fact, one can derive
rather explicit asymptotic expansions of integrals obtained by
integration along the fibers of a holomorphic function
(see [Bar82]).\qed

\section{\S2. $L^2$ extension theorem and inversion of adjunction}

Our starting point is the following special case of the fundamental
$L^2$ extension theorem due to Ohsawa-Takegoshi ([OhT87], [Ohs88], see
also [Man93]).

\claim 2.1. Theorem {\rm([OhT87], [Ohs88], [Man93])}|Let $\Omega\subset\bC^n$
be  a bounded pseudoconvex do\-main, and let $L$ be an affine linear subspace
of $\bC^n$ of codimension~$p\ge 1$ given by an orthonormal system $s$ of
affine linear equations $s_1=\cdots=s_p=0$. For every $\beta<p$, there
exists a constant $C_{\beta,n,\Omega}$ depending only on $\beta$, $n$
and the dia\-me\-ter of $\Omega$, satisfying the follo\-wing
property. For every $\varphi\in\cP(\Omega)$ and $f\in\cO(\Omega\cap L)$
with \hbox{$\int_{\Omega\cap L}|f|^2e^{-\varphi}dV_L<+\infty$}, there exists 
an extension $F\in\cO(\Omega)$ of $f$ such that
$$
\int_\Omega|F|^2|s|^{-2\beta}e^{-\varphi}dV_{\bC^n}\le C_{\beta,n,\Omega}
\int_{\Omega\cap L}|f|^2e^{-\varphi}dV_L,
$$
where $dV_{\bC^n}$ and $dV_L$ are the Lebesgue volume elements in $\bC^n$
and $L$ respectively.
\endclaim

In the sequel, we use in an essential way the fact that $\beta$
can be taken arbitrarily close to~$p$. It should be observed, however, that
the case $\beta=0$ is sufficient to imply the general case. In fact,
supposing $L=\{z_1=\cdots=z_p=0\}$, a substitution $(\varphi,\Omega)
\mapsto(\varphi_k,\Omega_k)$ with
$$
\eqalign{
&\varphi_k(z_1,\,\ldots\,,z_n)=\varphi(z_1^k,\,\ldots\,,z_p^k,z_{p+1},
\,\ldots\,,z_n),\cr
&\Omega_k=\big\{z\in\bC^n\,;\,(z_1^k,\,\ldots\,,z_p^k,z_{p+1},
\,\ldots\,,z_n)\in\Omega\big\}\cr}
$$
shows that the estimate with $\beta=0$ implies the estimate with
$\beta=p(1-1/k)$ (use the change of variable
$\zeta_1=z_1^k,\,\ldots\,,\zeta_p=z_p^k$, $\zeta_j=z_j$ for $j>p$,
together with the Jacobian formula
$$
dV(z)={\Const\over|\zeta_1|^{2(1-1/k)}\cdots|\zeta_p|^{2(1-1/k)}}
dV(\zeta),
$$
and take the ``trace'' of the solution $F_k$ on $\Omega_k$ to get the
solution $F$ on~$\Omega$). The $L^2$ extension theorem readily implies the
following important monotonicity result.

\claim 2.2. Proposition|Let $\varphi\in\cP(X)$ be a psh function on
a complex manifold~$X$, and let $Y\subset X$ be a complex submanifold such
that $\varphi_{|Y}\not\equiv-\infty$ on every connected component of~$Y$.
Then, if $K$ is a compact subset of $Y$, we have
$$
c_K(\varphi_{|Y})\le c_K(\varphi).
$$
$($Here, of course, $c_K(\varphi)$ is computed on $X$, i.e., by means of
neighborhoods of $K$ in~$X)$.
\endclaim

\proof. By Prop.~1.2, we may assume that $K=\{y\}$ is a single point
in~$Y$. Hence, after a change of coordinates, we can suppose that $X$
is an open set in~$\bC^n$ and that $Y$ is an affine linear subspace. Let
$c<c_y(\varphi_{|Y})$ be given. There is a small ball $B=B(y,r)$ such that
$\int_{B\cap Y}e^{-2c\varphi}dV_Y<+\infty$. By the $L^2$ extension theorem
applied with $\beta=0$, $\Omega=B$, $L=Y$ and $f(z)=1$, we can find a
holomorphic function $F$ on $B$ such that $F(z)=1$ on $B\cap Y$ and
$\int_B|F|^2e^{-2c\varphi}dV_B<+\infty$. As $F(y)=1$, we infer
$c_y(\varphi)\ge c$ and the conclusion follows. It should be observed
that an algebraic proof exists when $\varphi$ is of the form $\log|g|$,
$g\in\cO(X)\,$; however that proof is rather involved. This is already
a good indication of the considerable strength of the $L^2$ extension
theorem (which will be crucial in several respects in the sequel).\qed

We now show that the inequality given by Proposition~2.2 can somehow be
reversed (Theorem~2.5 below). For this, we need to restrict ourselves
to a class of psh functions which admit a ``sufficiently good local
behavior'' (such restrictions were already made in [Dem87], [Dem93a]
to accommodate similar difficulties).

\claim 2.3. Definition|Let $X$ be a complex manifold. 
We denote by $\cP_h(X)$ the class of all plurisubharmonic functions
$\varphi$ on $X$ such that $e^\varphi$ is locally H\"older continuous
on~$X$, namely such that for every compact set $K\subset X$ there are
constants $C=C_K\ge 0$, $\alpha=\alpha_K>0$ with
$$
|e^{\varphi(x)}-e^{\varphi(y)}|\le C\,d(x,y)^\alpha,\qquad\forall x,y\in K,
$$
where $d$ is some Riemannian metric on~$X$. We say for simplicity that
such a function is a H\"older psh function.
\endclaim

\claim 2.4. Example|{\rm We are mostly interested in the case of functions 
of the form
$$
\varphi=\max_j~\log\Big(\sum_k\prod_l|f_{j,k,l}|^{\alpha_{j,k,l}}\Big)
$$
with $f_{j,k,l}\in\cO(X)$ and $\alpha_{j,k,l}>0$. Such functions are easily
seen to be H\"older psh. Especially, if $D=\sum\alpha_jD_j$ is
an effective real divisor, the potential $\varphi=\sum\alpha_j\log|g_j|$
associated with $[D]$ is a H\"older psh function.}
\endclaim

\claim 2.5. Theorem|Let $H$ be a smooth hypersurface of $X$ and let $T$ be 
a closed positive current of type $(1,1)$ on $X$
such that its local potential functions $\varphi$ are H\"older psh functions
with $\varphi_{|H}\not\equiv-\infty$. We set in this case 
$($somewhat abusively$)$ $T_{|H}=dd^c\varphi_{|H}$.
Then for any compact set $K\subset H$, we have 
$$
c_K([H]+T)\ge 1 \quad \Leftrightarrow \quad c_K(T_{|H})\ge 1.
$$
\endclaim

In the algebraic setting (that is, when $T=[D]$ is defined by an
effective divisor $D=\sum\alpha_jD_j$), the above result is known as
``inversion of adjunction'', see Koll\'ar et al. [K\&al92], 17.7. One
says that the pair $(X,D)$ is lc (=~log canonical) if $c_K(D)\ge
1$ for every compact set $K\subset X$, i.e., if the product
$\prod|g_j|^{-2c\alpha_j}$ associated with the
generators $g_j$ of $\cO(-D_j)$ is locally $L^1$ for every $c<1$. 
The result can then be rephrased as
$$
\text{$(X,H+D)$ is lc}\quad \Leftrightarrow \quad \text{$(H,D_{|H})$ is lc}.
\leqno(2.5.1)
$$

\proof\ {\sl of Theorem 2.5}. Since the result is purely local, we may 
assume that $X=D(0,r)^n$ is a polydisk in $\bC^n$, that $H$ is the 
hyperplane $z_n=0$ and $K=\{0\}$. We must then prove the equivalence
$$
\eqalign{
\forall c<1,~\exists U\ni 0,~~&
\exp\big(-2c(\log|z_n|+\varphi(z))\big)\in L^1(U)\cr
&\Leftrightarrow \quad 
\forall c'<1,~\exists U'\ni 0,~~
\exp\big(-2c'\varphi(z',0)\big)\in L^1(U'),\cr}
$$
where $z=(z',z_n)\in\bC^n$ and $U$, $U'$ are neighborhoods of $0$ in $\bC^n$,
$\bC^{n-1}$ respectively.

First assume that $(|z_n|e^{\varphi(z)})^{-2c}\in L^1(U)$.
As $e^\varphi$ is H\"older continuous, we get 
$$
e^{2c\varphi(z)}\le(e^{\varphi(z',0)}+C_1|z_n|^{\alpha})^{2c}
\le C_2(e^{2c\varphi(z',0)}+|z_n|^{2c\alpha})
$$
on a neighborhood of $0$, for some constants $C_1$, $C_2$, $\alpha>0$. 
Therefore the function
$$
{1\over|z_n|^{2c}(|z_n|^{2c\alpha}+e^{2c\varphi(z',0)})}
\le C_2^{-1}(|z_n|e^{\varphi(z)})^{-2c}
$$
is in $L^1(U)$. Suppose that $U=U'\times D(0,r_n)$ is a small polydisk.
A partial integration with respect to $z_n$ on a family of disks 
$|z_n|<\rho(z')$ with $\rho(z')=\varepsilon\exp(\alpha^{-1}\varphi(z',0))$ 
(and $\varepsilon>0$ so small that $\rho(z')\le r_n$ for all $z'\in U'$)
shows that
$$
\int_U {dV(z)\over|z_n|^{2c}(|z_n|^{2c\alpha}+e^{2c\varphi(z',0)})}
\ge\text{Const}\int_{U'}{dV(z')\over e^{(2c-2(1-c)\alpha^{-1})\varphi(z',0)}}.
$$
Hence $\exp(-2c'\varphi(z',0))\in L^1(U')$ with $c'=c-(1-c)\alpha^{-1}$
arbitrarily close to~$1$. Conversely, if the latter condition holds,
we apply the Ohsawa-Takegoshi extension theorem to the function $f(z')=1$
on $L=H=\{z_n=0\}$, with the weight $\psi=2c'\varphi$ and $\beta=c'<1$.
Since $F(z',0)=1$, the $L^2$ condition implies the desired conclusion.\qed

\claim 2.6. Remark|{\rm As the final part of the proof shows, the implication
$$
c_K([H]+T)\ge 1 \quad \Leftarrow \quad c_K(T_{|H})\ge 1.
$$
is still true for an arbitrary (not necessarily H\"older) psh
function~$\varphi$. The implication $\Rightarrow$, however, is no longer
true. A simple counterexample is provided in dimension $2$ by $H=\{z_2=0\}$
and $T=dd^c\varphi$ with
$$
\varphi(z_1,z_2)=\max\big(\lambda\log|z_1|,{}-\sqrt{-\log|z_2|}\big),
\qquad \lambda>1
$$
on the unit bidisk $D(0,1)^2\subset\bC^2$. Then $c_0([H]+T)=c_0([H])=1$
but $c_0(T_{|H})=c_0(\lambda\log|z_1|)=1/\lambda$.}
\endclaim

\claim 2.7. Proposition|Let $X$, $Y$ be complex manifolds of 
respective dimensions $n$,~$m$, let $\cI\subset\cO_X$, $\cJ\subset\cO_Y$
be coherent ideals, and let $K\subset X$, $L\subset Y$ be compact
sets. Put $\cI\oplus\cJ:=\pr_1^\star\cI + \pr_2^\star\cJ\subset
\cO_{X\times Y}$. Then
$$
c_{K\times L}(\cI\oplus\cJ)=c_K(\cI)+c_L(\cJ).
$$
\endclaim

\proof. By Prop.~1.2, it is enough to show that
$c_{(x,y)}(\cI\oplus\cJ)=c_x(\cI)+c_y(\cJ)$ at every point $(x,y)\in
X\times Y$. Without loss of generality, we may assume that
\hbox{$X\subset\bC^n$}, \hbox{$Y\subset\bC^m$} are open sets and
$(x,y)=(0,0)$. Let $g=(g_1,\,\ldots\,,g_p)$, resp.\
$h=(h_1,\,\ldots\,,h_q)$, be systems of generators of $\cI$ (resp.\ $\cJ$)
on a neighborhood of~$0$. Set
$$
\varphi=\log\sum|g_j|,\qquad\psi=\log\sum|h_k|.
$$
Then $\cI\oplus\cJ$ is generated by the $p+q$-tuple of functions
$$
g\oplus h=(g_1(x),\,\ldots\,g_p(x),h_1(y),\,\ldots\,,h_q(y))
$$
and the corresponding psh function $\Phi(x,y)=\log\big(\sum|g_j(x)|+
\sum|h_k(y)|\big)$ has the same behavior along the poles as
$\Phi'(x,y)=\max(\varphi(x),\psi(y))$ (up to a term $O(1)\le\log 2$).
Now, for sufficiently small neighborhoods $U$, $V$ of $0$, we have
$$
\mu_{U\times V}\big(\big\{\max(\varphi(x),\psi(y))<\log r\big\}\big)
=\mu_U\big(\{\varphi<\log r\}\times\mu_V(\{\psi<\log r\}\big),
$$
hence Prop.~1.7$\,$(3) implies
$$
C_1r^{2(c+c')}\le
\mu_{U\times V}\big(\big\{\max(\varphi(x),\psi(y))<\log r\big\}\big)
\le C_2r^{2(c+c')}\,|\log r|^{n-1+m-1}
\leqno(2.7.1)
$$
with $c=c_0(\varphi)=c_0(\cI)$ and $c'=c_0(\psi)=c_0(\cJ)$. From this,
we infer 
$$
c_{(0,0)}(\cI\oplus\cJ)=c+c'=c_0(\cI)+c_0(\cJ).
\eqno\square
$$

\claim 2.8. Example|{\rm As $c_0(z_1^m)=1/m$, an application of
Proposition~2.7 to a quasi-homogeneous ideal $\cI=(z_1^{m_1}\ld z_p^{m_p})
\subset\cO_{\bC^n,0}$ yields the value
$$
c_0(\cI)={1\over m_1}+\cdots+{1\over m_p}.\eqno\square
$$
}
\endclaim

Using Proposition~2.7 and the monotonicity property, we can now prove
the fundamental subadditivity property of the singularity exponent.

\claim 2.9. Theorem|Let $f$, $g$ be holomorphic on a complex manifold~$X$.
Then, for every $x\in X$,
$$
c_x(f+g)\leq c_x(f)+c_x(g).
$$
More generally, if $\cI$ and $\cJ$ are coherent ideals, then
$$
c_x(\cI+\cJ)\leq c_x(\cI)+c_x(\cJ).
$$
\endclaim

\proof. Let $\Delta$ be the diagonal in $X\times X$. Then $\cI+\cJ$ can
be seen as the restriction of $\cI\oplus\cJ$ to~$\Delta$. Hence
Prop.~2.2 combined with 2.7 implies
$$
c_x(\cI+\cJ)=c_{(x,x)}((\cI\oplus\cJ)_{|\Delta})\le
c_{(x,x)}(\cI\oplus\cJ)=c_x(\cI)+c_x(\cJ).
$$
Since $(f+g)\subset(f)+(g)$, inequality~1.4$\,$(2) also shows that
$$
c_x(f+g)\le c_x((f)+(g))\le c_x(f)+c_x(g).\eqno\square
$$

\claim 2.10. Remark|{\rm If $f(x_1,\dots,x_n)$, resp.\ $g(y_1,\dots,y_n)$,
are holomorphic near $0\in\bC^n$, resp.\ $0\in \bC^m$, and such that
$f(0)=g(0)=0$, we have the equality
$$
c_0(f(x_1,\ldots,x_n)+g(y_1,\ldots,y_m))=\min\{1,c_0(f)+c_0(g)\}.
$$
This result is proved in [AGV84], vol.~II, sec.~13.3.5 in the case of
isolated singularities. Another proof, using the computation 
of $c_0$ via a resolution as in Prop.~1.7, is given in [Kol97]. It
can also be reduced to Proposition~2.7 through a log resolution of
either $f$ or~$g$.\qed}
\endclaim

\section{\S3. Semi-continuity of holomorphic singularity exponents}

We first give a new proof (in the spirit of this work) of the
semi-continuity theorem of Varchenko [Var82] concerning leading zeroes
of Bernstein-Sato polynomials attached to singularities of holomorphic
functions (see also Lichtin [Lin87]).

\claim 3.1. Theorem {\rm ([Var82])}|Let $X$ be a complex manifold
and $S$ a reduced complex space. Let $f(x,s)$ be a holomorphic function
on $X\times S$. Then for any $x_0\in X$, the function
$s\mapsto c_{x_0}(f_{|X\times \{s\}})$ is lower semi-continuous for the
holomorphic Zariski topology on~$S$. It even satisfies the 
following much stronger property: for any $s_0\in S$, one has
$$
c_{x_0}(f_{|X\times \{s\}})\ge c_{x_0}(f_{|X\times \{s_0\}})\leqno(3.1.1)
$$ 
on a holomorphic Zariski neighborhood of $s_0$ $($i.e.\ the complement in 
$S$ of an analytic subset of $S$ disjoint from~$s_0)$.
\endclaim

\proof. Observe that if $f_{|X\times\{s_0\}}$ is identically zero, then
$c_{x_0}(f_{|X\times \{s_0\}})=0$ and there is nothing to prove; thus we
only need to consider those $s$ such that $f_{|X\times\{s\}}\not\equiv 0$.
We may of course assume that $X=B$ is a ball in $\bC^n$ and~$x_0=0$. Let
$Y=B\times S$, $D=\div f$ and $\mu:\wt Y\to Y$ a log resolution of 
$(Y,D)$. After possibly shrinking $B$ a little bit, there is a Zariski 
dense open set $S_1\subset S$ such that if $s\in S_1$, the
corresponding fiber
$$
\mu_s:\wt Y_s\to B\times \{s\}
$$
is a log resolution of $(B,\div f_{|B\times \{s\}})$. Moreover, we may
assume that the numerical invariants  $a_i$, $b_i$ attached to
$\mu_s:\wt Y_s\to B$ as in Prop.~1.7 also do not depend on $s$.
In particular, by (1.7.2), $c_0(f_{|B\times \{s\}})$ is independent 
of $s\in S_1$.

By induction on the dimension of $S$, we obtain a stratification
$S=\bigcup S_i$ (where each $S_i$ is a Zariski dense open subset of a
closed complex subspace of $S$) such that $c_0(f_{|B\times \{s\}})$
only depends on the stratum containing~$s$. Thus (3.1.1) reduces to
semi-continuity with respect to the classical topology (considering a
$1$-dimensional base is enough, so we may assume the base to be
nonsingular as well). If we put $\varphi=\log|f|$, this is a special 
case of the following Lemma, which is essentially equivalent to the Main 
Theorem of [PS99]. Here, we would like to point out that this result
(which we knew as early as end of 1995) can be obtained as a direct
consequence of the Ohsawa-Takegoshi theorem [OhT87].

\claim 3.2. Lemma|Let $\Omega\subset\bC^n$ and $S\subset\bC^p$ be bounded
pseudoconvex open sets. Let $\varphi(x,s)$ be a H\"older psh function on 
$\Omega\times S$ and let $K\subset\Omega$ be a compact set. Then 
\smallskip
\item{\rm (1)} $s\mapsto c_K(\varphi(\bu,s))$ is lower semi-continuous for
the classical topology on~$S$.
\smallskip
\item{\rm (2)} If $s_0\in S$ and $c<c_K(\varphi(\bu,s_0))$, there exists
a neighborhood $U$ of $K$ and a uniform bound
$$
\int_U e^{-2c\varphi(x,s)}dV(x)\le M(c)
$$
for $s$ in a neighborhood of~$s_0$.
\vskip0pt
\endclaim

\proof. We use the $L^2$ extension theorem of [OhT87], following 
an idea of Angehrn-Siu [AnSi95]. However, the ``effective'' part (2)
requires additional considerations. Notice that it is enough to prove
(2), since (1) is a trivial consequence. By shrinking $\Omega$ and $S$,
we may suppose that $e^\varphi$ is H\"older continuous of exponent
$\alpha$ on the whole of $\Omega\times S$ and that
$$
\int_\Omega e^{-2c\varphi(x,s_0)}dV(x)<+\infty.
$$
Let $k$ be a positive integer. We set
$$
\psi_{k,s}(x,t)=2c\,\varphi(x,s+(kt)^k(s_0-s))
\qquad\text{on $\Omega\times D$,}
$$
where $D\subset\bC$ is the unit disk. Then $\psi$ is well defined on
$\Omega\times D$ if $s$ is close enough to~$s_0$. Since $\psi(x,1/k)
=\varphi(x,s_0)$, we obtain by Theorem~2.1 the existence of a holomorphic 
function $F_{k,s}(x,t)$ on $\Omega\times D$ such that $F_{k,s}(x,1/k)=1$
and
$$
\int_{\Omega\times D}|F_{k,s}(x,t)|^2e^{-\psi_{k,s}(x,t)}dV(x)dV(t)\le C_1
\leqno(3.2.3)
$$
with $C_1$ independent of~$k$,~$s$ for $|s-s_0|< \delta k^{-k}$. As 
$\psi_{k,s}$ admits a global
upper bound independent of~$k$,~$s$, the family $(F_{k,s})$ is a normal
family. It follows from the equality $F_{k,s}(x,1/k)=1$ that there is a
neighborhood $U$ of $K$ and a neighborhood $D(0,\varepsilon)$ of $0$
in $\bC$ such that $|F_{k,s}|\ge 1/2$ on $U\times D(0,\varepsilon)$ if
$k$ is large enough. A change of variable $t=k^{-1}\tau^{1/k}$ in
(3.2.3) then yields
$$
\int_{U\times D(0,(k\varepsilon)^k)}{e^{-2c\varphi(x,s+\tau(s_0-s))}
\over|\tau|^{2(1-1/k)}}\,dV(x)dV(\tau)\le 4k^4C_1.
$$
As in the proof of Theorem~2.5, we get by the H\"older continuity of 
$e^\varphi$ an upper bound
$$
e^{2c\varphi(x,s+\tau(s_0-s))}\le C_2(e^{2c\varphi(x,s)}+|\tau|^{2c\alpha})
$$
with a constant $C_2$ independent of~$s$. Hence, for $k\ge 1/\varepsilon$,
we find
$$
\int_{U\times D}
{1\over\big(e^{2c\varphi(x,s)}+|\tau|^{2c\alpha}\big)|\tau|^{2(1-1/k)}}\,
dV(x)dV(\tau)\le C_3(k).
$$
By restricting the integration to a family of disks 
$|\tau|<C_4e^{\alpha^{-1}\varphi(x,s)}$ (with $C_4$ so small that the 
radius is${}\le 1$), we infer
$$
\int_U
e^{-2(c-1/k\alpha)\varphi(x,s)}dV(x)\le C_5(k).
$$
Since $c-1/k\alpha$ can be taken arbitrarily close to $c_K(\varphi)$, this
concludes the proof.\qed

We can now prove the qualitative part of the semi-continuity theorem,
in the holomorphic case.

\claim 3.3. Theorem|Let $X$ be a complex manifold and $K\subset X$ a compact
subset. Then $f\mapsto c_K(f)$ is lower semi-continuous on $\cO(X)$ with
respect to the topology of uniform convergence on compact subsets. More
explicitly, for every nonzero holomorphic function $f$, for every compact
set $L$ containing $K$ in its interior and every $\varepsilon>0$, there is
a number $\delta=\delta(f,\epsilon,K,L)>0$ such that 
$$
\sup_L|g-f|<\delta \quad \Rightarrow \quad
c_K(g)\ge c_K(f)-\varepsilon.\leqno(3.3.1)
$$
\endclaim

\proof. As a first step we reduce (3.3.1) to the special case when $K$
is a single point. Assume that (3.3.1) fails. Then there is a sequence
of holomorphic functions $f_i\in\cO(X)$ converging uniformly to $f$ on
$L$, such that
$$
c_K(f_i)<c_K(f)-\varepsilon.
$$
By Prop.~1.2 we can choose for each $i$ a point $a_i\in K$ such that
$c_{a_i}(f_i)<c_K(f)-\varepsilon$. By passing to a subsequence we may
assume that the points $a_i$ converge to a point $a\in K$. Take a local
coordinate system on $X$ in a neighborhood of~$a$. Consider the
functions $F_i$ defined by
$$
F_i(x)=f_i(x+a_i-a)
$$
on a small coordinate ball $\ol B(a,r)\subset L^\circ$. These functions
are actually well defined for $i$ large enough (choose $\varepsilon$ so
that $\ol B(a,r+\varepsilon)\subset L$ and $i$ so large that
$|a_i-a|<\varepsilon$). Then $F_i$ converges to $f$ on $\ol B(a,r)$, but
$$
c_a(F_i)=c_{a_i}(f_i)<c_K(f)-\varepsilon\le c_a(f)-\varepsilon.
$$
Therefore, to get a contradiction, we only need proving Theorem~3.3 in
case $K=\{a\}$ is a single point. Again we can change notation and
assume that $X$ is the unit ball and that our point is the origin $0$.

In the second step we reduce the lower semi-continuity of $c_0(f)$ to
polynomials of bounded degree. For a given holomorphic function $f$ let
$P_k$ denote the degree${}\le k$ part of its Taylor series. The
subbaditivity property of Theorem~2.9 implies
$|c_0(f)-c_0(p_k)|\le c_0(f-p_k)$. As $|f(z)-p_k(z)|=O(|z|^{k+1})$, the
function $|f-p_k|^{-2c}$ is not integrable for $c\ge n/(k+1)$. From this,
it follows that $c_0(f-p_k)\le n/(k+1)$, hence
$$
|c_0(f)-c_0(p_k)|\le{n\over k+1}.\leqno(3.3.2)
$$
Now, if $(f_i)$ converges uniformly to $f$ on a given neighborhood
$U\subset\bC^n$ of~$0$, the degree${}\le k$ part $p_{i,k}$ converges to
$p_k$ in the finite dimension space $\bC[z_1\ld z_n]_k$ of polynomials
of total degree${}\le k$. Let us view polynomials
$$
P(z,s)=\sum_{|\alpha|\le k}s_\alpha z^\alpha\in \bC[z_1\ld z_n]_k
$$
as functions of their coefficients $s=(s_\alpha)$. By Theorem 3.1, we know
that the function $s\mapsto c_0(P(\bu,s))$ is lower semi-continuous.
Hence we get
$$
c_0(p_{i,k})>c_0(p_k)-{\varepsilon\over 2}
\quad\text{for $i>i(k,\varepsilon)$ large enough,}
$$
and thanks to (3.3.2) this implies
$$
c_0(f_i)>c_0(f)-{\varepsilon\over 2}-{2n\over k+1}>c_0(f)-\varepsilon
$$
by choosing $k\ge 4n/\varepsilon$.\qed

In fact, we would like to propose the following much stronger 
lower semi-continuity conjecture: 

\claim 3.4. Conjecture|Notation as in Theorem~3.3. For every nonzero
holomorphic function~$f$, there is a number $\delta=\delta(f,K,L)>0$
such that 
$$
\sup_L|g-f|<\delta \quad \Rightarrow \quad c_K(g)\ge c_K(f).
$$
\endclaim

\claim 3.5. Remark|{\rm There is an even more striking conjecture about
the numbers $c_K(f)$, namely, that the set 
$$
\cC=\{c_0(f)\vert f\in \cO_{\bC^n,0}\}\subset \bR
$$
satisfies the ascending chain condition (cf.\ [Sho92]; [K\&al92], 18.16):
any convergent increasing sequence in $\cC$ should be stationary. This 
conjecture and Theorem~3.3 together would imply the stronger form~3.4.
Notice on the other hand that there do exist non stationary decreasing
sequences in $\cC$ by~(1.4.8)\footnote{*}{\eightpoint
It has been recently observed by Phong and Sturm [PS00], in their study
of integrals of the form $\int|f|^{-s}$, that the ascending chain 
condition holds in complex dimension $2$. Algebraic geometers seem to
have been aware for some time of the corresponding algebraic
geometric statement.}.}
\endclaim

\section{\S4. Multiplier ideal sheaves and holomorphic approximations of 
psh singularities}

The most important concept relating psh functions to holomorphic objects
is the concept of {\sl multiplier ideal sheaf}, which was already
considered implicitly in the work of Bombieri [Bom70], Skoda [Sko72] and
Siu [Siu74]. The precise final formalization has been fixed
by Nadel [Nad89].

\claim 4.1. Theorem and definition {\rm ([Nad89,~90], see also [Dem89,~93a])}|
If $\varphi\in\cP(X)$ is a psh function on a complex manifold~$X$, the 
multiplier ideal sheaf $\cI(\varphi)\subset\cO_X$ is defined by
$$
\Gamma(U,\cI(\varphi))=
\big\{f\in\cO_X(U)\,;\,|f|^2e^{-2\varphi}\in L^1_\loc(U)\big\}
$$
for every open set $U\subset X$. Then $\cI(\varphi)$ is a coherent ideal
sheaf in $\cO_X$.
\endclaim

The proof that $\cI(\varphi)$ is coherent is a rather simple
consequence of H\"ormander's $L^2$ estimates, together with the strong
Noetherian property of coherent sheaves and the Krull lemma. When the
psh function $\varphi$ is defined from holomorphic functions as
in~2.4, it is easy to see that $\cI(\varphi)$ can be computed in a
purely algebraic way by means of log resolutions. The concept of
multiplier ideal sheaf plays a very important role in algebraic
geometry, e.g.\ in Nadel's version of the Kawamata-Viehweg vanishing
theorem or in Siu's proof [Siu93] of the big Matsusaka theorem.

We now recall the technique employed in [Dem92] and [Dem93b] to produce
effective bounds for the approximation of psh functions by logarithms
of holomorphic functions.  The same technique produces useful
comparison inequalities for the singularity exponents of a psh
function and its associated multiplier ideal sheaves.

\claim 4.2. Theorem|Let $\varphi$ be a plurisubharmonic function on a
bounded open set~$\Omega\subset\bC^n$. For every real number $m>0$, let
$\cH_{m\varphi}(\Omega)$ be the Hilbert space of holomorphic functions $f$
on $\Omega$ such that $\int_\Omega|f|^2e^{-2m\varphi}dV<+\infty$ and
let $\psi_m={1\over 2m}\log\sum|g_{m,k}|^2$ where $(g_{m,k})$
is an orthonormal basis of~$\cH_{m\varphi}(\Omega)$. Then\/$:$
\smallskip
\item{\rm (1)} There are constants $C_1,C_2>0$ independent of $m$ and
$\varphi$ such that
$$
\varphi(z)-{C_1\over m}\le
\psi_m(z)\le\sup_{|\zeta-z|<r}\varphi(\zeta)+{1\over m}\log{C_2\over r^n}
$$
for every $z\in\Omega$ and $r<d(z,\partial\Omega)$. In particular,
$\psi_m$ converges to $\varphi$ pointwise and in $L^1_{\rm loc}$ topology
on~$\Omega$ when $m\to+\infty$ and
\smallskip
\item{\rm (2)} The Lelong numbers of $\varphi$ and $\psi_m$ are related by
$$
\nu(\varphi,z)-{n\over m}\le\nu(\psi_m,z)\le
\nu(\varphi,z)\quad\text{ for every $z\in\Omega$.}
$$
\item{\rm (3)} For every compact set $K\subset\Omega$, the Arnold 
multiplicity of $\varphi$,~$\psi_m$ and of the multiplier ideal sheaves
$\cI(m\varphi)$ are related by
$$
\lambda_K(\varphi)-{1\over m}\le\lambda_K(\psi_m)=
{1\over m}\lambda_K(\cI(m\varphi))\le \lambda_K(\varphi).
$$
\vskip0pt
\endclaim


\proof. (1) Note that $\sum|g_{m,k}(z)|^2$ is the square of the norm of 
the evaluation linear form $f\mapsto f(z)$ on $\cH_{m\varphi}(\Omega)$. As
$\varphi$ is locally bounded above, the $L^2$ topology is actually
stronger than the topology of uniform convergence on compact subsets
of~$\Omega$.  It follows that the series $\sum|g_{m,k}|^2$ converges
uniformly on $\Omega$ and that its sum is real analytic. Moreover we have
$$
\psi_m(z)=\sup_{f\in B(1)}{1\over m}\log|f(z)|
$$
where $B(1)$ is the unit ball of $\cH_{m\varphi}(\Omega)$. For
$r<d(z,\partial\Omega)$, the mean value inequality applied to the psh
function $|f|^2$ implies
$$
\eqalign{
|f(z)|^2&\le{1\over\pi^nr^{2n}/n!}\int_{|\zeta-z|<r}
|f(\zeta)|^2d\lambda(\zeta)\cr
&\le{1\over\pi^nr^{2n}/n!}\exp\Big(2m\sup_{|\zeta-z|<r}\varphi(\zeta)\Big)
\int_\Omega|f|^2e^{-2m\varphi}d\lambda.\cr}
$$
If we take the supremum over all $f\in B(1)$ we get
$$
\psi_m(z)\le\sup_{|\zeta-z|<r}\varphi(\zeta)+{1\over 2m}
\log{1\over\pi^nr^{2n}/n!}
$$
and the right hand inequality in (1) is proved. Conversely, the
Ohsawa-Takegoshi extension theorem applied to the $0$-dimensional
subvariety $\{z\}\subset\Omega$ shows that for any $a\in\bC$ there is
a holomorphic function $f$ on $\Omega$ such that $f(z)=a$ and
$$
\int_\Omega|f|^2e^{-2m\varphi}d\lambda\le C_3|a|^2e^{-2m\varphi(z)},
$$
where $C_3$ only depends on $n$ and~$\diam\Omega$. We fix $a$
such that the right hand side is~$1$. This gives the left hand
inequality
$$
\psi_m(z)\ge{1\over m}\log|a|=\varphi(z)-{\log C_3\over 2m}.
\leqno(4.2.4)
$$
(2) The above inequality (4.2.4) implies $\nu(\psi_m,z)\le\nu(\varphi,z)$. 
In the opposite direction, we find
$$
\sup_{|x-z|<r}\psi_m(x)\le\sup_{|\zeta-z|<2r}\varphi(\zeta)+
{1\over m}\log{C_2\over r^n}.
$$
Divide by $\log r$ and take the limit as $r$ tends to~$0$. The quotient
by $\log r$ of the supremum of a psh function over $B(x,r)$ tends to the
Lelong number at~$x$. Thus we obtain
$$
\nu(\psi_m,x)\ge\nu(\varphi,x)-{n\over m}.
$$
(3) Inequality (4.2.4) already yields
$\lambda_K(\psi_m)\le\lambda_K(\varphi)$.  Moreover, the multiplier
ideal sheaf $\cI(m\varphi)$ is generated by the sections in
$\cH_{m\varphi}(\Omega)$ (as follows from the proof that
$\cI(m\varphi)$ is coherent), and by the strong Noetherian property,
it is generated by finitely many functions $(g_{m,k})_{0\le k\le k_0(m)}$ 
on every relatively compact open set $\Omega'\compact\Omega$.
It~follows that we have a lower bound of the form
$$
\psi_m(z)-C_4\le
{1\over 2m}\log\sum_{0\le k\le k_0(m)}|g_{m,k}|^2\le\psi_m(z)\quad
\text{on $\Omega'$.}\leqno(4.2.5)
$$
By choosing $\Omega'\supset K$, we infer $\lambda_K(\psi_m)={1\over
  m}\, \lambda_K(\cI(m\varphi)$. If $\lambda>\lambda_K(\psi_m)$, i.e.,
$1/m\lambda<c_K(\cI(m\varphi))$, and if $U\subset\Omega'$ is a
sufficiently small open neighborhood of $K$, the H\"older inequality
for the conjugate exponents $p=1+m\lambda$ and $q=1+(m\lambda)^{-1}$
yields
$$
\leqalignno{
\int_U e^{-2mp^{-1}\varphi}dV
&=\int_U\Big(\sum_{0\le k\le k_0(m)}|g_{m,k}|^2e^{-2m\varphi}\Big)^{1/p}
  \Big(\sum_{0\le k\le k_0(m)}|g_{m,k}|^2\Big)^{-1/qm\lambda}dV\cr
&\le (k_0(m)+1)^{1/p}
\left(\int_U \Big(\sum_{0\le k\le k_0(m)}|g_{m,k}|^2\Big)^{-1/m\lambda}dV
\right)^{1/q}<+\infty.&(4.2.6)
\cr}
$$
$\big($The estimate in the last line uses the fact that
$$
\int_U|g_{m,k}|^2e^{-2m\varphi}dV\le 
\int_\Omega|g_{m,k}|^2e^{-2m\varphi}dV=1.~\big)
$$
This implies $c_K(\varphi)\ge mp^{-1}$, i.e., $\lambda_K(\varphi)\le p/m=
\lambda+1/m$. As $\lambda>\lambda_K(\psi_m)$ was arbitrary, we get
$\lambda_K(\varphi)\le\lambda_K(\psi_m)+1/m$ and (3) follows.\qed

The ``approximation theorem'' 4.2 allows to extend some results proved for
holomorphic functions to the case of psh functions. For instance, we have:

\claim 4.3. Proposition|Let $\varphi\in\cP(X)$, $\psi\in\cP(Y)$ be psh 
functions on complex manifolds $X$,~$Y$, and let $K\subset X$, $L\subset Y$ 
be compact subsets. Then\/$:$
\smallskip
\item{\rm(1)} For all positive real numbers $c',\,c''$ with 
$c'>c_K(\varphi)>c''$ $($if any$)$ and every sufficiently small 
neighborhood $U$ of $K$, there is an estimate
$$
C_1r^{2c'}\le\mu_U(\{\varphi<\log r\})\le C_2r^{2c''},
\qquad \forall r<r_0
$$
for some $r_0>0$ and $C_1=C_1(c')$, $C_2=C_2(c'')$.
\smallskip
\item{\rm(2)} $c_{K\times L}\big(\max(\varphi(x),\psi(y)\big)=
c_K(\varphi)+c_L(\psi)$.
\smallskip
\item{\rm(3)} If $X=Y$, then $c_x(\max(\varphi,\psi))
\le c_x(\varphi)+c_x(\psi)$ for all $x\in X$.\vskip0pt
\endclaim

\proof. (1) The upper estimate is clear, since
$$
r^{-2c''}\mu_U(\{\varphi<\log r\})\le\int_Ue^{-2c''\varphi}dV<+\infty
$$
for $U\subset K$ sufficiently small. In the other direction, we have an 
estimate
$$
\mu_U(\{\psi_m<\log r\})\ge C_{1,m}r^{2c_K(\psi_m)}
$$
by Proposition 1.7$\,$(3) and (4.2.5). As $\varphi\le\psi_m+C_{2,m}$ for some
constant $C_{2,m}>0$, we get
$$
\{\varphi<\log r\}\supset\{\psi_m<\log r-C_{2,m}\},
$$ 
and as $c_K(\psi_m)$ converges to $c_K(\varphi)$ by 4.2$\,$(3), the
lower estimate of $\mu_U(\{\varphi<\log r\})$ follows.

\noindent(2), (3) can be derived from (1) exactly as for the holomorphic
case in Prop.~2.7 and Theorem~2.9. It should be observed that 4.3$\,$(1)
expresses a highly non trivial ``regularity property'' of the growth of
volumes $\mu_U(\{\varphi<\log r\})$ when $\varphi$ is a psh function
(when $\varphi$ is an arbitrary measurable function, $v(r)=
\mu_U(\{\varphi<\log r\})$ is just an arbitrary increasing function with
$\lim_{r\to 0}v(r)=0$).\qed

\claim 4.4. Remark|{\rm In contrast with the holomorphic case
1.7$\,$(3), the upper estimate $\mu_U(\{\varphi<\log r\})\le
C_2r^{2c''}$ does not hold with $c''=c_K(\varphi)$, when $\varphi$ is
an arbitrary psh function. A simple example is given by
$\varphi(z)=\chi\circ\log|z|$ where $\chi:\bR\to\bR$ is a convex
increasing function such that $\chi(t)\sim t$ as $t\to-\infty$, but
$e^{\chi(r)}\not\sim r$ as $r\to 0$, e.g.\ such that
$\chi(t)=t-\log|t|)$ when $t<0$. On the other hand, the
lower estimate $\mu_U(\{\varphi<\log r\})\ge C_1r^{2c'}$ seems to be
still true with $c'=c_K(\varphi)$, although we cannot prove~it.}
\endclaim

\section{\S5. Semi-continuity of psh singularity exponents}

We are now in a position to prove our main semi-continuity theorem.

\noindent {\bf 5.1.} {\sl Proof of Theorem 0.2}. Let
$\Omega\subset\bC^n$ be a bounded pseudoconvex open set and let
$\varphi_j\in\cP(\Omega)$ be a sequence of psh functions converging to
a limit $\varphi\in\cP(\Omega)$ in the weak topology of distributions.
In fact, this already implies that $\varphi_j\to\varphi$ almost everywhere 
and in $L^1_\loc$ topology; to see this, we observe that the coefficients
of $T_j:=dd^c\varphi_j$ are measures converging to those of $T=dd^c\varphi$
in the weak topology of measures; moreover $\varphi_j$ and $\varphi$ can
be recovered from $T_j$ and $T$ by an integral formula involving the
Green kernel; we then use the well known fact that integral operators
involving a $L^1$ kernel define continuous (and even compact) operators
from the space of positive measures equipped with the weak topology,
towards the space of $L^1$ functions with the strong $L^1$ topology. 

By the process described in
Theorem~4.2, we get for each $m\in\bN^\star$ a Hilbert orthonormal basis
$(g_{j,m,k})_{k\in\bN}$ of $\cH_{m\varphi_j}(\Omega)$, such that
$$
\varphi_j(z)-{C_1\over m}\le{1\over 2m}\log\sum_{k\in\bN}|g_{j,m,k}|^2
\le\sup_{|\zeta-z|<r}\varphi_j(\zeta)+{1\over m}
\log{C_2\over r^n}\leqno(5.1.1)
$$
for every $z\in\Omega$ and $r<d(z,\partial\Omega)$. 
\RGBColor{1 0 0}{[Additional explanations added 
on November~13, 2013, following a discussion with JingZhou Sun].}
\RGBColor{0 0 1}{We need here a special choice of the orthonormal 
basis $(g_{j,m,k})_{k\in\bN}$ to ensure convergence as $j\to+\infty$,
since especially the sequence could be randomly permuted with respect
to the index~$k$. To this end, select a sequence of points 
$(a_\ell)_{\ell\in\bN}$
that is dense in $\Omega$, such that all $a_\ell$ are in the complement of
the negligible set $V(\cI(\varphi))\cup\bigcup_jV(\cI(\varphi_j))\subset
\Omega$, i.e.\ such that all functions $e^{-\varphi}$ and $e^{-\varphi_j}$ are 
locally integrable in a neighborhood of~$a_\ell$. Let $\cS_{j,m,k}\subset
\cH_{m\varphi_j}(\Omega)$ be the (closed) codimension $k$ subspace of functions 
$f\in \cH_{m\varphi_j}(\Omega)$ such that $f(a_\ell)=0$ for $0\leq\ell<k$
(the assertion about codimension follows from a standard interpolation argument
based on solving a $\dbar$-equation with H\"ormander's $L^2$ estimates;
however, we do not really need this argument, it would be enough to observe
that $\codim\cS_{j,m,k}\leq k$ is finite, and for this, any dense sequence
$(a_\ell)$ would work).
Since $\bigcap_{k\in\bN}\cS_{j,m,k}=\{0\}$ by density of the sequence $(a_\ell)$,
we can select the orthonormal basis $(g_{j,m,k})_{k\in\bN}$ of
$\cH_{m\varphi_j}(\Omega)$ so that
$g_{j,m,k}\in \cS_{j,m,k}\cap \cS_{j,m,k+1}^\perp$. By (5.1.1),
the Cauchy-Schwarz inequality and the existence of a local uniform 
upper bound for the $\varphi_j$, the Bergman kernel function
$$
\beta_{j,m}(z,w)=\sum_{k\in\bN}g_{j,m,k}(z)\overline{g_{j,m,k}(w)}
$$
is uniformly bounded from above on every compact subset of $\Omega\times\Omega$.
As $\beta_{j,m}(z,w)$ is holomorphic in $(z,\overline{w})$, one can extract 
a subsequence $p\mapsto j_p$ such that $\smash{\widetilde\beta_{m}(z,w)}
:=\lim_{p\to +\infty}\beta_{j_p,m}(z,w)$ locally 
uniformly on $\Omega\times\Omega$, with a limit that is again
holomorphic in $(z,\overline{w})$. By restricting to the diagonal
$z=w$, one sees that the individual sequences
$(g_{j,m,k})_{j\in\bN}$ are also locally uniformly bounded from above.
After extracting a diagonal subsequence still denoted
$p\mapsto j_p$, we may assume that all $g_{j,m,k}$ converge to a
limit $g_{m,k}\in\cO(\Omega)$ when $j=j_p\to+\infty$. An application of
the discrete Fatou lemma implies
$$
\widetilde\beta'_m(z,z):=
\sum_{k\in\bN}|g_{m,k}(z)|^2\le \widetilde\beta_m(z,z)=\lim_{p\to +\infty}
\sum_{k\in\bN}|g_{j_p,m,k}(z)|^2.
$$
However, we have equality for $z=a_\ell$, since the value
$\beta_{j,m}(a_\ell,a_\ell)$ is already achieved by truncating the series
as a finite sum $\smash{\sum_{0\le k\le\ell}}$ 
(by construction $g_{j,m,k}(a_\ell)=0$ for $k>\ell$). Since 
$\smash{\widetilde\beta_m}$ and $\smash{\widetilde\beta'_m}$ are real
analytic and coincide on the dense sequence $(a_\ell)$, we must have 
$\smash{\widetilde\beta_m(z)=\widetilde\beta'_m(z)}$ everywhere.
Thanks to (5.1.1) we find in the limit
$$
\varphi(z)-{C_1\over m}\le{1\over 2m}\log\sum_{k\in\bN}|g_{m,k}(z)|^2
\le\sup_{|\zeta-z|<r}\varphi(\zeta)+{1\over m}
\log{C_2\over r^n}.
$$
(Notice that if two psh functions $u$, $v$ satisfy
$u\leq v$ almost everywhere, then $u\le v$ everywhere).}

Fix a compact set $K\subset\Omega$ and a relatively compact open 
subset $\Omega'\compact\Omega$
containing~$K$. By the strong Noetherian property already used for (4.2.5), 
there exist an integer $k_0(m)$ and a constant $C_4(m)>0$ such that
$$
\varphi(z)-C_4(m)\le{1\over 2m}\log\sum_{0\le k\le k_0(m)}|g_{m,k}(z)|^2\quad
\text{on $\Omega'$.}
$$
Now, for $c<c_K(\varphi)$, there is a neighborhood $U$ of $K$ on which
$$
\int_U\Big(\sum_{0\le k\le k_0(m)}|g_{m,k}|^2\Big)^{-c/m}dV\le
e^{2cC_4(m)}\int_Ue^{-2c\varphi}dV<+\infty.
$$
Take (without loss of generality) $m\ge 2c_K(\varphi)$. Then $c/m<1/2$ and
Formula~1.5.1 shows that there is a linear combination
$\sum_{0\le k\le k_0(m)}\alpha_{m,k}g_{m,k}$ with $\alpha=(\alpha_{m,k})$
in the unit sphere of $\bC^{k_0(m)+1}$, such that
$$
\int_U\Big|\sum_{0\le k\le k_0(m)}\alpha_{m,k}g_{m,k}\Big|^{-2c/m}dV\le
C_5(m)\int_Ue^{-2c\varphi}dV<+\infty,
$$
where $C_5(m)$ is a constant depending possibly on~$m$. By construction,
$$
f_{j,m}=\sum_{0\le k\le k_0(m)}\alpha_{k,m}g_{j,m,k}
$$
is an element of the unit sphere in $\cH_{m\varphi_j}(\Omega)$, and
$f_{j,m}$ converges uniformly on $\Omega$ to $f_m=\sum\alpha_{m,k}g_{m,k}$
such that $\int_U|f_m|^{-2c/m}dV<+\infty$. By Lemma~5.2 below, for
any $c'<c$ and $K\subset U'\compact U$, we have a uniform bound
$\int_{U'}|f_{j,m}|^{-2c'/m}dV\le C_6(m)$ for $j\ge j_0$ large enough.
Since $\int_\Omega|f_{j,m}|^2e^{-2m\varphi_j}dV=1$, the H\"older inequality
for conjugate exponents $p=1+m/c'$, $q=1+c'/m$ yields
$$
\eqalign{
\int_{U'}e^{-2mc'/(m+c')\varphi_j}dV
&=\int_{U'}\big(|f_{j,m}|^2e^{-2m\varphi_j}\big)^{c'/(m+c')}
  |f_{j,m}|^{-2c'/(m+c')}dV\cr
&\le\Big(\int_{U'}|f_{j,m}|^{-2c'/m}dV\Big)^{m/(m+c')}\le C_7(m)\cr}
$$
for $j\ge j_0$. Since $c$, $c'$ are arbitrary with $c'<c<c_K(\varphi)$,
the exponent $mc'/(m+c')$ can be taken to approach $c$ as closely as we
want as $m$ gets large. Hence $c_K(\varphi_j)>c_K(\varphi)-\varepsilon$
for $j\ge j_0(\varepsilon)$ large enough. Moreover, by what we have seen 
above, if $c<c_K(\varphi)$ is fixed and $0<\delta<c_K(\varphi)/c-1$, there
exists $j_1(\delta)$ such that the sequence $(e^{-2c\varphi_j})_{j\ge
j_1(\delta)}$ is contained in a bounded set of $L^{1+\delta}(U)$,
where $U$ is a small neighborhood of~$K$. Therefore
$$
\int_{U\cap\{e^{-2c\varphi_j}>M\}}e^{-2c\varphi_j}dV\le C_8M^{-\delta}
$$
for $j\ge j_1(\delta)$, with a constant $C_8$ independent of~$j$.
Since $e^{-2c\varphi_j}$ converges pointwise to $e^{-2c\varphi}$ on $\Omega$,
an elementary argument based on Lebesgue's bounded convergence theorem 
shows that $e^{-2c\varphi_j}$ converges to $e^{-2c\varphi}$ in~$L^1(U)$.\qed

To complete the proof, we need only proving the following effective
estimate for holomorphic functions, which is a special case
of part (3) in the Main Theorem.

\claim 5.2. Lemma|Let $\Omega\subset\bC^n$ be a bounded pseudoconvex open
set, and let $f_i\in\cO(\Omega)$ be a sequence of holomorphic 
functions converging uniformly to $f\in\cO(\Omega)$ on every compact 
subset. Fix a compact set $K\subset\Omega$ and $c<c_K(f)$. Then there 
is a neighborhood $U$ of $K$ and a uniform bound $C>0$ such that
$$
\int_U |f_i|^{-2c}dV\le C
$$ 
for $i\ge i_0$ sufficiently large.
\endclaim

\proof. We already know by Theorem 3.3 that $\int_U |f_i|^{-2c}dV<+\infty$ 
for $U$ small enough and $i$~large. Unfortunately, the proof given in 
Theorem~3.3 is not effective because it depends (through the use
of Hironaka's theorem in the proof of estimates 1.7$\,$(3) and (2.7.1)) 
on the use of a sequence of log resolutions on which we have absolutely 
no control. We must in fact produce an effective version of 
inequality~(3.3.2).

The result of Lemma~5.2 is clearly local. Fix a point $x_0\in K$
(which we assume to be $0$ for simplicity), real numbers $c'$, 
$c''$ with $c<c''<c'<c_K(f)\le c_0(f)$ and an integer $k$ so large that 
$$
c<c''-{n\over k+1}<c''<c'<c_0(f)-{n\over k+1}.
$$
Let $p_k$ be the truncation at order $k$ of the Taylor series of $f$ at the
origin. As $c_0(p_k)\ge c_0(f)-{n\over k+1}>c'$ by (3.3.2), there is a small 
ball $B'=B(0,r_0')$ such that
$$
\int_{B'}|p_k|^{-2c'}dV<+\infty.
$$
Since the truncations $p_{i,k}$ of $f_{i,k}$ converge uniformly to $p_k$
on $\bC^n$ as $i\to+\infty$, Lemma~3.2 applied to the universal family
of polynomials $P(z,s)=\sum_{|\alpha|\le k}s_\alpha z^\alpha$ shows that 
for any ball $B''\compact B'$, there is a constant $M\ge 0$ and an 
integer $i_0$ such that
$$
\int_{B''}|p_{i,k}|^{-2c''}dV\le M\qquad\text{for $i\ge i_0$}.
$$
Let us write $p_{i,k}=f_i-g_{i,k}$ where $g_{i,k}$ consists of the sum of
terms of degree${}>k$ in the Taylor expansion of $f_i$ at the origin.
By the Ohsawa-Takegoshi theorem applied with the weight function
$\psi(x,y)=2c\log|f_i(x)-g_{i,k}(y)|$ on $B''\times B''$ and
$L={}$diagonal of $\bC^n\times\bC^n$, there is a holomorphic function
$F_i$ on $B''\times B''$ such that $F_i(x,x)=1$ and
$$
\int_{B''\times B''}|F_i(x,y)|^2|f_i(x)-g_{i,k}(y)|^{-2c''}dV(x)\,dV(y)\le C_1
$$
with a constant $C_1$ independent of $i$. The above $L^2$ estimate shows
that $(F_i)$ is bounded in $L^2$ norm on $B''\times B''$. Hence, there
is a small ball $B=B(0,r_0)\compact B''$ such that $|F_i(x,y)|\ge 1/2$ 
on $B\times B$ for all~$i\ge i_0$, and
$$
\int_{B\times B}|f_i(x)-g_{i,k}(y)|^{-2c''}dV(x)\,dV(y)\le 4C_1.\leqno(5.2.1)
$$
Moreover, we have a uniform estimate $|g_{i,k}(y)|\le C_2|y|^{k+1}$ on~$B$
with a constant $C_2$ independent of~$i$. By integrating (5.2.1) with 
respect to $y$ on the family of balls $|y|<(|f_i(x)|/2C_2)^{1/(k+1)}$,
we find an estimate
$$
\int_B|f_i(x)|^{2n/(k+1)-2c''}dV(x)\le C_3.\leqno(5.2.2)
$$
As $c''-n/(k+1)>c$, this is the desired estimate. It is interesting to 
observe that the proof of the Main Theorem can now be made entirely
independent of Hironaka's desingularization theorem. In fact, the only
point where we used it is in the inequality $c_0(p_k)\ge c_0(f)
-{n\over k+1}$, which we derived from Proposition~2.7. The latter 
inequality can however be derived directly from the Ohsawa-Takegoshi 
theorem through estimates for 
$\int_{B\times B}|p_k(x)+g_k(y)|^{-2c}dV(x)\,dV(y)$.\qed

\claim 5.3. Remark|{\rm It follows from the proof of Prop.~1.7 that the
set of positive exponents $c$ such that $|f|^{-2c}$ is summable on a
neighborhood of a compact set $K$ is always an {\sl open interval}, namely
$]0,c_K(f)[$. We conjecture that the same property holds true more generally
for an arbitrary psh function $\varphi$ (``openness conjecture'');
the openness conjecture is indeed true in dimension $1$, since we have
the well known necessary and sufficient criterion
$$
e^{-2\varphi}\in L^1_\loc(V(x_0))~~
\Leftrightarrow~~ \nu(\varphi,x_0)<1
$$
(as follows e.g.\ from [Sko72]). By using the Main Theorem, the openness
conjecture would imply the following stronger statement:}
\endclaim

\claim 5.4. Strong openness conjecture|Let $U'\compact U\compact X$
be relatively compact open sets in a complex manifold~$X$. Let $\varphi$
be a psh function on~$X$ such that $\int_Ue^{-\varphi}dV<+\infty$. Then
there exists $\varepsilon=\varepsilon(\varphi,U,U')$ such that for
every $\psi$ psh on~$X$
$$
\Vert\psi-\varphi\Vert_{L^1(U)}<\varepsilon \quad\Rightarrow\quad
\int_{U'}e^{-\psi}dV<+\infty.
$$
In other words, the integrability of $e^{-\varphi}$ near a given compact
set $K$ should be an open property for the $L^1_\loc$ $(=\,$ weak$)$
topology on $\cP(X)$.
\endclaim

The main theorem only yields the weaker conclusion
$$
\int_{U'}e^{-(1-\delta)\psi}dV<+\infty\quad\hbox{for}\quad
\Vert\psi-\varphi\Vert_{L^1(U)}<\varepsilon=\varepsilon(\varphi,U,U',\delta).
$$

\section{\S6. Existence of K\"ahler-Einstein metrics on Fano orbifolds}

An {\sl orbifold} is a complex variety $X$ possessing only quotient 
singularities, namely, every point $x_0\in X$ has a neighborhood $U$
isomorphic to a quotient $\Omega/\Phi$ where $\Phi=\Phi_{x_0}$ is a
finite group acting holomorphically on a smooth open set
$\Omega\subset\bC^n$. Such an action can always be linearized, so we may
assume that $\Phi$ is a finite subgroup of $\GL_n(\bC)$ and
$\Omega$ a $\Phi$-invariant neighborhood of~$0$ (with $x_0$ being the
image of $0$). We may also assume that the elements of $G$
distinct from  identity have a set of fixed points of codimension~${}\ge 2$
(otherwise, the subgroup generated by these is a normal subgroup $N$
of $\Phi$, $\Omega/N$ is again smooth, and $\Omega/\Phi=
(\Omega/N)/(\Phi/N)$). The structure sheaf $\cO_X$ (resp.\ the $m$-fold
canonical sheaf $K_X^{\otimes m}$) is then defined locally as the direct
image by $\pi:\Omega\to U\simeq\Omega/\Phi$ of the subsheaf of
$\Phi$-invariant sections of the corresponding sheaf on $\Omega\,:$
$$
\Gamma(V,\cO_X)=\Gamma(\pi^{-1}(V),\cO_\Omega)^\Phi,\qquad
\Gamma(V,K_X^{\otimes m})=\Gamma(\pi^{-1}(V),K_\Omega^{\otimes m})^\Phi,
$$
for all open subsets $V\subset U$. There is always an integer $m_0$ (e.g.\
$m_0=\#\,\Phi$) such that $K_\Omega^{\otimes m_0}$ has $\Phi$-invariant
local generating sections, and then clearly $K_X^{\otimes m}$ is an 
invertible $\cO_X$-module whenever $m$ is divisible by the lowest 
common multiple $\mu$ of the integers $m_0$ occurring in the various 
quotients.
Similarly, one can define on $U$ (and thus on $X$) the concepts of
K\"ahler metrics, Ricci curvature form, etc, by looking at corresponding
$\Phi$-invariant objects on~$\Omega$. We say that a compact
orbifold $X$ is a {\sl Fano orbifold} if $K_X^{-\mu}$ is ample, which
is the same as requiring that $K_X^{-\mu}$ admits a smooth hermitian
metric with positive definite curvature. In that case, we define the
curvature of $K_X^{-1}$ to be $1/\mu$ times the curvature of $K_X^{-\mu}$.
The integral of a differential form on $X$ (say defined at least 
on $X_\reg$) is always computed upstairs, i.e.\ 
$\int_{\Omega/\Phi}\alpha={1\over\#\Phi}\int_\Omega\pi^\star\alpha$.

\claim 6.1. Definition|A compact orbifold $X$ is said to be
K\"ahler-Einstein if it possesses a K\"ahler form $\omega
={i\over 2\pi}\sum\omega_{jk}dz_j\wedge d\ol z_k$ satisfying
the Einstein condition
$$
\Ricci(\omega)=\lambda\omega
$$
for some real constant $\lambda$, where where $\Ricci(\omega)$ is the
closed $(1,1)$-form defined in every coordinate patch by
$\Ricci(\omega)=-{i\over 2\pi}\ddbar\log\det(\omega_{jk})$.
\endclaim

Since $\Ricci(\omega)$ is the curvature form of $K_X^{-1}=\det T_X$ 
equipped with the metric $\det\omega$, a necessary condition for the
existence of a K\"ahler-Einstein metric with constant $\lambda>0$ is
that $K_X^{-1}$ is ample, i.e., that $X$ is Fano. On the other hand, it
is well known that not all Fano orbifolds are K\"ahler-Einstein, even
when they are smooth; further necessary conditions are required, e.g.\
that the group of automorphisms $\Aut(X)^\circ$ is reductive ([Mat57],
[Lic57]), and that the {\sl Futaki invariants} vanish [Fut83]; for
instance $\bP^2$ blown up in $2$ points has a non reductive group of
automorphisms and therefore is not K\"ahler-Einstein.

It is usually much harder to prove that a concretely given Fano orbifold
is K\"ahler-Einstein. Siu [Siu87,~88], and slightly later Tian [Tia87]
and Nadel [Nad89,~90], gave nice sufficient conditions ensuring the
existence of a K\"ahler-Einstein metric; these conditions always involve
the existence of a sufficiently big group of automorphisms. Our goal here
is to reprove Nadel's main result in a more direct and conceptual way.

\claim 6.2. Technical setting|{\rm We first briefly recall the main
technical tools and notation involved (see e.g.\ [Siu87] for more
details). The anticanonical line bundle $K_X^{-1}$ is assumed to
be ample. Therefore it admits a smooth hermitian metric $h_0$
whose $(1,1)$-curvature form $\theta_0={i\over 2\pi}D_{h_0}^2$ is
positive definite. Since $\theta_0\in c_1(X)$, the Aubin-Calabi-Yau
theorem shows that there exists a K\"ahler metric $\omega_0\in c_1(X)$
such that $\Ricci(\omega_0)=\theta_0$. [The Aubin-Calabi-Yau is still valid
in the orbifold case, because the proof depends only on local regularity 
arguments which can be recovered by passing to a finite cover, and global
integral estimates which still make sense by the remark preceding 
Def.~6.1]. Since both $\theta_0$ and $\omega_0$ are in $c_1(X)$, we have
$$
\omega_0=\theta_0+{i\over 2\pi}\ddbar f\qquad
\hbox{for some $f\in C^\infty(X)$.}\leqno(6.2.1)
$$
We look for a new K\"ahler form $\omega=\omega_0+{i\over 2\pi}\ddbar\varphi$
in the same K\"ahler class as $\omega_0$, such that $\Ricci\omega=\omega$.
Since $\Ricci(\omega_0)=\theta_0$, this is equivalent to
$$
-{i\over 2\pi}\ddbar\log(\det\omega)=\omega=\theta_0+{i\over 2\pi}
\ddbar(\varphi+f)=-{i\over 2\pi}\ddbar\log(\det\omega_0)
+{i\over 2\pi}\ddbar(\varphi+f),
$$
that is,
$$
\ddbar\Big(\log{\det\omega\over\det\omega_0}+\varphi+f\Big)=0,
$$
which in its turn is equivalent to the Monge-Amp\`ere equation
$$
\log{(\omega_0+{i\over 2\pi}\ddbar\varphi)^n\over \omega_0^n}
+\varphi+f+C=0
\leqno(6.2.2)
$$
where $C$ is a constant. Here, one can normalize $\varphi$ so that
$\varphi$ is orthogonal to the \hbox{$1$-dimensional} space of constant
functions in $L^2(X,\omega_0)$, i.e., $\int_X\varphi\,\omega_0^n=0$. The
usual technique employed to solve (6.2.2) is the so-called ``continuity
method''. The continuity method amounts to introducing an extra parameter
$t\in[0,1]$ and looking for a solution $(\varphi_t,C_t)$ of the equation
$$
\log{(\omega_0+{i\over 2\pi}\ddbar\varphi_t)^n\over \omega_0^n}
+t(\varphi_t+f)+C_t=0,\qquad \int_X\varphi_t\omega_0^n=0
\leqno(6.2.3)
$$
as $t$ varies from $0$ to~$1$. Clearly $\varphi_0=0$, $C_0=0$ is a
solution for $t=0$ and $(\varphi,C)=(\varphi_1,C_1)$ provides a
solution of our initial equation (6.2.2). Moreover, the linearization
of the (nonlinear) elliptic differential operator occuring in (6.2.3)
is the operator
$$
(\psi,C)\longmapsto{1\over 2\pi}\Delta_{\omega_t}\psi +t\psi+C
\leqno(6.2.4)
$$
where $\omega_t$ is the K\"ahler metric $\omega_t=\omega_0+{i\over2\pi}
\ddbar\varphi_t$ and $\Delta_{\omega_t}$ is the associated Laplace 
operator (with negative eigenvalues). The equation (6.2.3) is easily seen 
to be equivalent to
$$
\Ricci(\omega_t)=t\omega_t+(1-t)\theta_0.
$$ 
From this we infer $\Ricci(\omega_t)>t\omega_t$ for all $t<1$, and it 
then follows from
the Bochner-Kodaira-Nakano identity that all nonzero eigenvalues
of $-{1\over 2\pi}\Delta_{\omega_t}$ are ${}>t$ (this is clear directly
for $-{1\over 2\pi}\Delta_{\omega_t}$ acting on $(0,1)$-forms, and one
uses the fact that $\dbar$ maps the $\lambda$-eigenspace
$E^{p,q}(\lambda)$ of $-{1\over 2\pi}\Delta_{\omega_t}$ in bidegree
$(p,q)$ into the corresponding eigenspace $E^{p,q+1}(\lambda)$).
Then, thanks to Schauder's estimates, (6.2.4) induces an isomorphism
$\cC^{s+2}_\perp(X)\oplus\bR\to\cC^s(X)$ where $s\in\bR_+\ssm\bN$
and $\cC^s(X)$ (resp.\ $\cC^s_\perp(X)$) is the space of real
functions (resp.\ real functions orthogonal to constants) of
class $\cC^s$ on~$X$. Let $\cT\subset[0,1]$ be the set of parameters
$t$ for which (6.2.3) has a smooth solution. By elliptic regularity
for (nonlinear) PDE equations, the existence of a smooth solution is
equivalent to the existence of a solution in $\cC^s(X)$ for some $s>2$.
It then follows by a standard application of the implicit function
theorem that $\cT\cap[0,1[$ is an open subset of the interval~$[0,1[$.}
\endclaim

\claim 6.3. Sufficient condition for closedness|{\rm In order to obtain
a solution for all times $t\in[0,1]$, one still has to prove that $\cT$
is {\it closed}. By the well-known theory of complex Monge-Amp\`ere
equations ([Aub78], [Yau78]), a sufficient condition for closedness
is the existence of a uniform a priori $\cC^0$-estimate
$\Vert\wt\varphi_t\Vert_{\cC^0}\le \Const$ for the family of functions
$\wt\varphi_t=t\varphi_t+C_t$, $t\in\cT$, occuring in the right hand
side of~(6.2.3). A first observation is that
$$
\sup_X\varphi_t\le \Const,\qquad\text{hence}~~\sup_X\wt\varphi_t
\le C_t+\Const,
\leqno(6.3.1)
$$
as follows from the conditions $\int_X\varphi_t\omega_0^n=0$ and
${i\over 2\pi}\ddbar\varphi_t\ge-\omega_0$, by simple considerations of
potential theory. On the other hand, by [Siu88, Prop.~2.1] or
[Tia87, Prop.~2.3], we have the Harnack-type inequality
$$
\sup_X(-\wt\varphi_t)\le (n+\varepsilon)\sup_X\wt\varphi_t+A_\varepsilon,
\leqno(6.3.2)
$$
where $\varepsilon>0$ and $A_\varepsilon$ is a constant depending only on
$\varepsilon$. Hence
$\sup_X(-\wt\varphi_t)\le(n+\varepsilon)C_t+A'_\varepsilon$ and we thus
only need controlling the constants~$C_t$ from above. Now, equation
(6.2.3) implies
$$
\int_X\omega_0^n=\int_X\Big(\omega_0+{i\over 2\pi}
\ddbar\varphi_t\Big)^n=\int_X e^{-\wt\varphi_t-tf}\omega_0^n.
$$
For $\gamma\in{}]0,1[$, we easily infer from this and (6.3.2) that
$$
\eqalign{\int_X\omega_0^n
&\le\Const\,\exp\big((1-\gamma)\sup_X(-\wt\varphi_t)\big)
    \int_X e^{-\gamma\wt\varphi_t}\omega_0^n\cr
&\le\Const_\varepsilon\,e^{(1-\gamma)(n+\varepsilon)C_t}
    \int_X e^{-\gamma\wt\varphi_t}\omega_0^n\cr
&\le\Const_\varepsilon\,e^{-(\gamma-(1-\gamma)(n+\varepsilon))C_t}
    \int_X e^{-\gamma t\varphi_t}\omega_0^n.\cr}
$$
If $\gamma\in{}]{n\over n+1},1[$ and $\varepsilon$ is small enough,
we conclude that $C_t$ admits an upper bound of the form
$$
C_t\le B'_\gamma\log\int_Xe^{-\gamma t\varphi_t}\omega_0^n+B''_\gamma
$$
where $B'_\gamma$ and $B''_\gamma$ depend only on~$\gamma$. Hence
closedness of $\cT$ is equivalent to the uniform boundedness of the
integrals
$$
\int_Xe^{-\gamma t\varphi_t}\omega_0^n,\qquad t\in\cT,
\leqno(6.3.3)
$$
for any choice of $\gamma\in{}]{n\over n+1},1[$.}
\endclaim

This yields the following basic existence criterion due to
Nadel [Nad89,~90].

\claim 6.4. Existence criterion for K\"ahler-Einstein metrics|Let $X$ be
a Fano orbifold of dimension~$n$. Let $G$ be a compact subgroup of the
group of complex automorphisms of~$X$. Then $X$ admits a $G$-invariant
K\"ahler-Einstein metric, unless $K_X^{-1}$ possesses a $G$-invariant
singular hermitian metric $h=h_0e^{-\varphi}$ $(h_0$~being a
smooth $G$-invariant metric and $\varphi$ a $G$-invariant function in
$L^1_\loc(X))$, such that the following properties occur.
\smallskip
\item{\rm(1)} $h$ has a semipositive curvature current
$$
\Theta_h=-{i\over 2\pi}\ddbar\log h=
\Theta_{h_0}+{i\over 2\pi}\ddbar\varphi\ge 0.
$$
\item{\rm(2)} For every $\gamma\in{}]{n\over n+1},1[$, the multiplier
ideal sheaf $\cI(\gamma\varphi)$ is nontrivial,  $($i.e.\
$0\ne\cI(\gamma\varphi)\ne\cO_X)$.
\vskip0pt
\endclaim

According to the general philosophy of orbifolds, the orbifold concept
of a multiplier ideal sheaf $\cI(\gamma\varphi)$ is that the ideal
sheaf is to be computed upstairs on a smooth local cover and take the
direct image of the subsheaf of invariant functions by the local
isotropy subgroup; this ideal coincides with the multiplier ideal
sheaf computed downstairs only if we take downstairs the volume form
which is the push forward of an invariant volume form upstairs (which
is in general definitely larger than the volume form induced by a
local smooth embedding of the orbifold).

\proof. Let us start with a $G$-invariant K\"ahler metric
$\omega_0={i\over 2\pi}\ddbar\log h_0^{-1}$, where $h_0$ and $\omega_0$
have the same meaning as in 6.2; indeed, if $h_0$ is not $G$-invariant,
we can average it by using the $G$-action, that is, we define a new
metric $(h_0^G)^{-1}$ on $K_X$ by putting
$$
(h_0^G)^{-1}=\int_{g\in G}g^\star h_0^{-1}d\mu(g),
$$
and we again have $\omega_0^G:={i\over 2\pi}\ddbar\log (h_0^G)^{-1}>0$.
Now, all $\varphi_t$ can be taken to be $G$-invariant. If the continuity
process ceases to produce a solution $\varphi_t$ at $t=t_0\in[0,1]$
(thus, if $t_0\in\ol\cT\ssm\cT$), there exists a sequence
$t_\nu\in\cT$ converging to $t_0$ and (6.3.3) implies
$\lim_{\nu\to+\infty}\int_Xe^{-\gamma t_\nu\varphi_{t_\nu}}
\omega_0^n=+\infty$ for every $\gamma\in{}]{n\over n+1},1[$.
As the space of closed positive currents contained in a given
cohomology class is compact for the weak topology, one can extract a
subsequence $\Theta_{(p)}=\omega_0+{i\over 2\pi}\ddbar\varphi_{t_{\nu(p)}}$
converging weakly to a limit $\Theta=\omega_0+{i\over 2\pi}\ddbar\varphi
\ge 0$. The potential $\varphi$ can be recovered from $\Tr\Theta$ by means
of the Green kernel, and therefore, by the well-known properties of the
Green kernel, we infer that $\varphi_{t_{\nu(p)}}$ converges to $\varphi$
in~$L^1(X)$. The semicontinuity theorem in its effective version 0.2.2
shows that
$$
\int_Xe^{-\gamma t_0\varphi}\omega_0^n=+\infty\qquad
\text{for all $\gamma\in{}]{n\over n+1},1[$,}
$$
and therefore $\int_Xe^{-\gamma\varphi}\omega_0^n=+\infty$
for all $\gamma\in{}]{n\over n+1},1[$. From this we conclude that
$\cI(\gamma\varphi)\ne \cO_X$. The fact that $\cI(\gamma\varphi)\ne 0$
is clear since $\varphi\not\equiv-\infty$.\qed

Before going further, we need Nadel's vanishing theorem (a generalized
version of the well-known Kawamata-Viehweg vanishing theorem. It is
known to be a rather simple consequence of H\"ormander's $L^2$
estimates, see e.g.\ [Dem89], [Nad89], [Dem93b] or [Dem94]).

\claim 6.5. Nadel vanishing theorem|Let $(X,\omega)$ be a K\"ahler
orbifold and let $L$ be a holomorphic orbifold line bundle over $X$ 
equipped with a singular hermitian metric $h$ of weight~$\varphi$ with
respect to a smooth metric $h_0$ $($i.e.\ $h=h_0e^{-\varphi})$.
Assume that the curvature form $\Theta_h(L)={i\over 2\pi}D_h^2$ is
positive definite in the sense of currents, i.e.\
$\Theta_h(L)\ge\varepsilon\omega$ for some $\varepsilon>0$. If 
$K_X\otimes L$ is an invertible sheaf on $X$, we have
$$
H^q\big(X,K_X\otimes L\otimes\cI(\varphi)\big)=0\qquad
\text{for all $q\ge 1$.}
$$
{\rm Recall that an ``orbifold line bundle'' $L$, is a rank $1$ sheaf 
which is locally an invariant direct image of an invertible sheaf on 
$\Omega$ by the local quotient maps $\Omega\to\Omega/\Phi\,$; $L$ itself 
need not be invertible; similarly, $\otimes$ is meant to be the orbifold
tensor product, i.e., we take the tensor product upstairs on $\Omega$
and take the direct image of the subsheaf of invariants. The proof is
obtained by the standard $L^2$ estimates applied on $X_\reg$ with respect
to an orbifold K\"ahler metric on $X$. It is crucial that $K_X\otimes L$ 
be invertible on $X$, otherwise the set of holomorphic sections of 
$K_X\otimes L$ satisfying the $L^2$ estimate with respect to the weight
$e^{-\varphi}$ might differ from the orbifold tensor product
$K_X\otimes L\otimes\cI(\varphi)$ [and also, that tensor product might 
be equal to $K_X\otimes L$ even though $\cI(\varphi)$ is non trivial].}
\endclaim

\claim 6.6. Corollary|Let $X$, $G$, $h$ and $\varphi$ be as in
Criterion~6.4. Then, for all $\gamma\in{}]{n\over n+1},1[$,
\smallskip
\item{\rm(1)} the multiplier ideal sheaf $\cI(\gamma\varphi)$ satisfies
$$
H^q(X,\cI(\gamma\varphi))=0\qquad\text{for all $q\ge 1$.}
$$
\smallskip
\item{\rm(2)} the associated subscheme $V_\gamma$ of structure
sheaf $\cO_{V_\gamma}=\cO_X/\cI(\gamma\varphi)$ is nonempty, distinct
from $X$, $G$-invariant and satisfies
$$
H^q(V_\gamma,\cO_{V_\gamma})
=\cases{
\bC&for $q=0$,\cr
0  &for $q\ge 1$.\cr}
$$
\vskip0pt
\endclaim

\proof. Apply Nadel's vanishing theorem to $L=K_X^{-1}$ equipped
with the singular hermitian metric $h_\gamma=h_0e^{-\gamma\varphi}$.
Then $\Theta_{h_\gamma}=\gamma\Theta_h+(1-\gamma)\Theta_{h_0}\ge
(1-\gamma)\omega_0>0$, and (1) follows. Finally, since $X$ is Fano,
we get
$$
H^q(X,\cO_X)=0\qquad\text{for all $q\ge 1$},
$$
by Kodaira vanishing for $L=K_X^{-1}$. The exact sequence
$$
0\to \cI(\gamma\varphi)\to\cO_X\to\cO_{V_\gamma}\to 0
$$
immediately implies (2).\qed

The strategy employed by Nadel [Nad90] to construct K\"ahler-Einstein
metrics is to rule out the existence of any $G$-invariant subscheme
with the properties described in 6.6$\,$(2). Of course, this is easier
to achieve if $G$ is large. One uses the following observations
(assuming that the closedness property fails, so that all subschemes
$V_\gamma$ are nontrivial).

\claim 6.7. Proposition|All subschemes $V_\gamma$ are connected.
Therefore, if $G$ has no fixed points, $V_\gamma$ cannot be
$0$-dimensional.
\endclaim

\proof. The connectedness of $V_\gamma$ is a straightforward consequence
of the equality \hbox{$H^0(V_\gamma,\cO_{V_\gamma})=\bC$}.\qed

\claim 6.8. Proposition|If $V_\gamma$ contains irreducible components
$Z_j$ of codimension $1$, then the corresponding divisor $Z=\sum m_jZ_j$
satisfies the numerical inequality $[Z]\le\gamma[K_X^{-1}]$ in the
sense that $\gamma[K_X^{-1}]-[Z]$ can be represented by a closed
positive current. In particular, one always has the inequality
$$
(-K_X)^{n-1}\cdot Z\le\gamma(-K_X)^n.
$$
If $K_X^{-1}$ generates the group $W(X)$ of Weil divisors of $X$ modulo
numerical equivalence, then $V_\gamma$ must have codimension $\ge 2$.
\endclaim

\noindent In the smooth case we have of course $W(X)=\Pic(X)$, but in general
$\Pic(X)$ is a subgroup of finite index in~$W(X)$.

\proof. Consider the closed positive $(1,1)$ current
$\Theta_h=\omega_0+{i\over 2\pi}\ddbar\varphi$ (which belongs to the
first Chern class $c_1(K_X^{-1})$), and let
$$
\Theta_h=\sum\lambda_j[Z_j]+R,\qquad\lambda_j>0,~~R\ge 0
$$
be the Siu decomposition of $\Theta_h$ (namely, the $[Z_j]\,$'s are
currents of integration over irreducible divisors and $R$ is a closed
$(1,1)$-current which has nonzero Lelong numbers only in
codimension~$2$). It is then easy to see that
the subscheme $V_\gamma$ defined by $\cI(\gamma\varphi)$ precisely has
$[Z]=\sum\lfloor\gamma\lambda_j\rfloor\,[Z_j]$ as its $1$-codimensional
part (here, $\lfloor~~\rfloor$ denotes the integral part). Hence
$\gamma \Theta_h-[Z]\ge 0$ as asserted. If $K_X^{-1}$ generates $\Pic(X)$,
this implies $Z=0$, since there cannot exist any nonzero effective
integral divisor numerically smaller than $[K_X^{-1}]$.\qed

When $\dim X=3$, $G$ has no fixed points and $K_X^{-1}$ generates
$W(X)$, we are only left with the case $V_\gamma$ is of pure
dimension~$1$. This case can sometimes be ruled out by observing that 
certain groups cannot act effectively on the curve $V_\gamma$
(As $H^1(V_\gamma,\cO_{V_\gamma})=0$, $V_\gamma$ is a tree of rational
curves; see Nadel [Nad90, Th.~4.1, 4.2 and Cor.~4.1]).

Further a priori inequalities can be derived for certain components
of the multiplier ideal subschemes~$V_\gamma$. Especially, for
components of codimension~$2$, we have the following simple bound,
based on a use of a self-intersection inequality for the current
$\Theta=\omega_0+{i\over 2\pi}\ddbar\varphi$.

\claim 6.9. Proposition|Assume that $W(X)$ is generated by $K_X^{-1}$
and that $a$ is a nonnegative number such that the orbifold vector bundle
$T_X\otimes\cO(-aK_X)$ is numerically effective. Then the codimension $2$
components $Z_j$ of $V_\gamma$ satisfy the inequality
$$
\sum {1\over \delta_j}\nu_j(\nu_j-1)(-K_X)^{n-2}\cdot Z_j\le(1+a)(-K_X)^n
$$
where $\nu_j\ge 1/\gamma$ is the generic Lelong number of $\Theta=\omega_0
+{i\over 2\pi}\ddbar\varphi$ along~$Z_j$, and $\delta_j$ is the order
of the local isotropy group of the orbifold at a generic point in~$Z_j$.
Especially, if $\gamma$ is taken
to be sufficiently close to ${n\over n+1}$, we have
$$
\sum(-K_X)^{n-2}\cdot Z_j\le{n^2\over n+1}\delta(1+a)(-K_X)^n
$$
where $\delta$ is the maximum of the the orders of the isotropy groups.
\endclaim

\proof. Since $V_\gamma$ is of codimension $2$ for $\gamma$ arbitrarily
close to~$1$, the generic Lelong number of $\varphi$ must be${}\le 1$
along all components of codimension $1$ in the Lelong sublevel sets
$E_c(\varphi)=E_c(\Theta)$ [again, Lelong numbers and Lelong sublevel sets
are to be interpreted upstairs, on a smooth finite cover]. If a 
codimension $2$ component $Z_j$ occurs
in $\cI(\gamma\varphi)$, the generic Lelong number $\gamma\nu_j$ of
$\gamma\varphi$ along that component must be $\ge 1$, hence
$\nu_j\ge 1/\gamma$. We now apply the regularization theorem for closed
$(1,1)$-currents ([Dem92], Main Theorem). For every
$c>1$ we obtain a current $\Theta_{h,c}$ cohomologous to $\Theta$
(hence in the class $c_1(K_X^{-1})$), which is smooth on
$X\ssm E_c(\Theta)$, thus smooth except on an analytic set of
codimension${}\ge 2$, such that $\Theta_{h,c}\ge -(ca+\varepsilon)
\omega_0$ and such that the Lelong numbers of $\Theta_{h,c}$ are
shifted by $c$, i.e.\ $\nu_x(\Theta_{h,c})=(\nu_x(\Theta)-c)_+$. The
intersection product
$\Theta\wedge (\Theta_c+(ca+\varepsilon)\omega_0)$ is well defined, belongs
to the cohomology class $(1+ca+\varepsilon)(-K_X)^2$ and is larger than
$\sum{1\over \delta_j}\nu_j(\nu_j-c)[Z_j]$ as a current. Hence, by taking 
the intersection with the class $(-K_X)^{n-2}$ we get
$$
\sum {1\over \delta_j}\nu_j(\nu_j-c)(-K_X)^{n-2}\cdot 
Z_j\le(1+ca+\varepsilon)(-K_X)^n.
$$
[The extra factor $1/\delta_j$ occurs because we have to divide by $\delta_j$ 
to convert an integral on a finite cover $\Omega$ to an integral on
the quotient $\Omega/\Phi$].
As $c$ tends to $1+0$ and $\varepsilon$ tends to $0+$, we get the desired
inequality. The last observation comes from the fact that $\cI(V_\gamma)$
must be constant on some interval $]{n\over n+1},{n\over n+1}+\delta[$,
by the Noetherian property of coherent sheaves.\qed
\bigskip

\claim 6.10. Example|{\rm Let $\bP_a=\bP^3(a_0,a_1,a_2,a_3)$ be the weighted
projective $3$-space with weights $a_0\le a_1\le a_2\le a_3$ such that the
 components $a_i$ are relatively prime $3$ by~$3$. It is equipped with
an orbifold line bundle $\cO_X(1)$ which, in general, is not locally free.
Let $t=a_0+a_1+a_2+a_3$ and
$$
X=\big\{P(x_0,x_1,x_2,x_3)=0\big\}
$$
be a generic surface of weighted degree $d$ in $\bP_a$.
It is known (see Fletcher [Fle89]) that $X$ has an orbifold structure 
(i.e., is quasi-smooth in the terminology of Dolgachev [Dol82]),
if and only if the following conditions are satisfied:
\smallskip
\item{(i)} For all $j$, there exists a monomial $x_j^mx_{k(j)}$ of degree 
$d\,$;
\smallskip
\item{(ii)} For all distinct $j$, $k$, either there exists a monomial
$x_j^mx_k^p$ of degree $d$, or there exist monomials
$x_j^{m_1}x_k^{p_1}x_{\ell_1}$, $x_j^{m_2}x_k^{p_2}x_{\ell_2}$ of degree $d$
with $\ell_1\ne\ell_2\,$;
\smallskip
\item{(iii)} For all $j$, there exists a monomial of degree $d$ which does not
involve $x_j$.
\smallskip
\noindent
 Moreover, $-K_X=\cO_X(t-d)$ (and  hence $(-K_X)^2=d(t-d)^2/(a_0a_1a_2a_3)$)
if and only if the following condition also holds:
 \smallskip
\item{(iv)} For every $j,k$ such that $a_j$ and $a_k$ are not relatively prime,
 there exists a monomial $x_j^mx_k^p$ of degree $d$.
\smallskip
\noindent
We would like to use the conditions of Propositions~6.8 and 6.9
to show that $X$ carries a K\"ahler-Einstein metric.

Proposition~6.8 clearly applies if we can prove that
$(-K_X)\cdot Z>{2\over 3} (-K_X)^2$ for every effective curve on $X$.
This is not a priori trivial in the examples below since the Picard 
numbers will always be bigger than $1$. Using the torus action, every curve 
on a weighted projective space can be 
degenerated to a sum of lines of the form $(x_i=x_j=0)$.
Thus $(-K_X)\cdot Z$ is bounded from below by $(t-d)/(a_2a_3)$.
Thus $(-K_X)\cdot Z>{2\over 3} (-K_X)^2$ holds if
$$
{t-d\over {a_2a_3}}> {2\over 3} {{d(t-d)^2}\over {a_0a_1a_2a_3}},
\qquad\hbox{i.e.}\quad
a_0a_1>{2\over 3}d(t-d).
$$
In the examples we give at the end, which all concern the case 
$d=t-1$,  this is always satisfied.

In order to apply Proposition~6.9, we need to determine $T_X$.
We have exact sequences
$$
\eqalign{
&0\to\cO_{\bP_a}\to\bigoplus\cO_{\bP_a}(a_i)\to T_{\bP_a}\to 0,\cr
&0\to T_X\to T_{\bP_a|X}\to \cO_X(d)\to 0,\cr
}
$$
and we get from there a surjective arrow
$$
\bigoplus\cO_X(a_i)\to \cO_X(d)
$$
given explicitly by the matrix $(\partial P/\partial x_i)$.
{}From the above exact sequences, we find a sequence of surjective arrows
$$
\bigoplus_{i<j}\cO_X(a_i+a_j)\to
\cO_X(\Lambda^2T_{\bP_a|X})\to T_X\otimes\cO_X(d).
$$
(Of course, formally speaking, we are dealing with orbifold vector bundles,
which can be considered as locally free sheaves only when we pass to
a finite Galois cover). 
Moreover,
$$
\bigoplus_{i\neq k\neq j}\cO_X(a_i+a_j)\to
 T_X\otimes\cO_X(d)
$$
is surjective over the open set where $x_k\neq 0$. 
This proves that, as an orbifold vector bundle, $T_X\otimes\cO_X(d-a_0-a_2)$
is nef if the line $(x_0=x_1=0)$ is not contained in~$X$.

The maximal order $\delta$ of the isotropy groups is less than 
$a_3$ --~which is indeed
the maximum for $\bP_a$ itself~-- resp.\ $a_2$ if $a_3$ divides $d$, since 
in that case a generic surface of degree $d$ does not pass through the 
point $[0:0:0:1]$. 
 This shows that we can take $a=(d-a_0-a_2)/(t-d)$ in Proposition~6.9,
and as the $Z_j$ are points and $n=2$, we find the condition
$$
1\le{4\over 3}a_3\Big(1+{d-a_0-a_2\over t-d}\Big)
{d(t-d)^2\over a_0a_1a_2a_3},
$$
with the initial $a_3$ being replaced by $a_2$ if $a_3$ divides $d$. We thus
compute the ratio
$$
\eqalign{
\rho_a&={4\over 3}{d(t-d)(t-a_0-a_2)\over a_0a_1a_2}
\quad\hbox{if $a_3\mathop{\not|}d$},\cr 
\cr
\rho_a&={4\over 3}{d(t-d)(t-a_0-a_2)\over a_0a_1a_3}\quad
\hbox{if $a_3|d$},\cr
}
$$
and when $\rho_a<1$ we can conclude that the Del Pezzo surface is 
K\"ahler-Einstein. Clearly, this is easier to reach when $t-d$ is
small, and we concentrated ourselves on the case $d=t-1$. 
It is then easy to check that $\rho_a$ is never less than
$1$ when $a_0=a_1=1$. On the other hand, a computer check seems to
indicate that there is only a finite list of weights with
$a_0>2$ satisfying the Fletcher conditions, which all 
satisfy $a_0\le 14$\footnote{*}{Added after proof:
this has actually been shown to be true in [JK01].}. 
Among these, 2 cases lead to $\rho_a<1$, namely
$$
\kern-10pt
\matrix{
&a=(11,49,69,128),\hfill&d=256,\hfill&\kern-5pt\rho_a\simeq 0.875696,&
x_0^{17}x_2+x_0x_1^5+x_1x_2^3+x_3^2=0,\hfill\cr
&a=(13,35,81,128),\hfill&d=256,\hfill&\kern-5pt\rho_a\simeq 0.955311,&
x_0^{17}x_1+x_0x_2^3+x_1^5x_2+x_3^2=0.\hfill\cr
}
$$
It turns out that there are no other monomials of degree $d$ than those 
occurring in the above equations. As a result, the above K\"ahler-Einstein
Del Pezzo surfaces are {\it rigid} as weighted hypersurfaces.

There are several ways to improve the estimates. For instance,
$$
T_X\otimes\cO_X(d-a_1-a_2)
$$
is nef except possibly along the irreducible components of the curve
$(x_0=0)\subset X$.  The restriction of the tangent bundle to these
curves can be computed by hand.  This improvement is sufficient to
conclude that Propositions~6.8 and 6.9 also apply in one further
case:
$$
\kern-10pt
\matrix{
&a=(9,15,17,20),\hfill&d=60,\hfill&
x_0^{5}x_1+x_0x_2^3+x_1^4+x_3^3=0.\hfill\cr
}
$$
This is again a rigid weighted hypersurface. We would like to thank
P.~Boyer and K.~Galicki for pointing out a numerical error which had
been committed in an earlier version of this work, where a further
(incorrect) example $a=(11,29,39,49)$, $d=127$ was claimed. In [BG00],
it is shown that the three above examples lead to the construction of
non regular Sasakian-Einstein $5$-manifolds.
}
\endclaim

\bigskip\null
\section{References}
\medskip

{\eightpoint
  
\bibitem[AnSi95]&U.\ Angehrn {\rm and} Y.-T.\ Siu:& Effective Freeness
  and Point Separation for Adjoint Bundles;& Invent.\ Math.\ 122
  (1995) 291--308&

\bibitem[AnV65]&A.\ Andreotti {\rm and} E.\ Vesentini:& Carleman
  estimates for the Laplace-Beltrami equation in complex manifolds;&
  Publ.\ Math.\ I.H.E.S.\ {\bf 25} (1965), 81--130&
  
\bibitem[ArGV85]&V.I.\ Arnold, S.M.\ Gusein-Zade {\rm and} A.N.\ 
  Varchenko;& Singularities of Differentiable Maps;& Progress in
  Math., Birkh\"auser (1985)&
  
\bibitem[Aub78]&T.\ Aubin:& Equations du type Monge-Amp\`ere sur les
  vari\'et\'es k\"ahl\'eriennes compactes;& C.R.\ Acad.\ Sci.\ Paris
  Ser.\ A {\bf 283} (1976), 119--121$\,$; Bull.\ Sci.\ Math.\ {\bf
  102} (1978), 63--95&

\bibitem[Bar82]&D.\ Barlet:& D\'eveloppements asymptotiques des fonctions
  obtenues par int\'egration sur les fibres;& Invent.\ Math.\ {\bf 68}
  (1982), 129--174&

\bibitem[Bom70]&E.\ Bombieri:& Algebraic values of meromorphic maps;&
  Invent.\ Math. {\bf 10} (1970), 267--287 and Addendum,
  Invent.\ Math.\ {\bf 11} (1970), 163--166&

\bibitem[BG00]&C.\ Boyer, K.\ Galicki:& New Sasakian-Einstein $5$-manifolds
  as links of isolated hypersurface singularities;& 
  Manuscript, February 2000&

\bibitem[Dem87]&J.-P.\ Demailly:& Nombres de Lelong g\'en\'eralis\'es,
  th\'eor\`emes d'int\'egralit\'e et d'analy\-ti\-cit\'e;& Acta Math.\ 
  {\bf 159} (1987) 153--169&

\bibitem[Dem89]&J.-P.\ Demailly:& Transcendental proof of a generalized
  Ka\-wa\-mata-Viehweg vanishing theorem;& C.\ R.\ Acad.\ Sci.\ Paris
  S\'er.\ I Math.\ {\bf 309} (1989) 123--126~ and~ 
  Proceedings of the Conference ``Geometrical and algebraical aspects
  in several complex variables'' held at Cetraro (Italy), C.A.~Berenstein
  and D.C.~Struppa eds, EditEl, June (1989), 81--94&
  
\bibitem[Dem90]&J.-P.\ Demailly:& Singular hermitian metrics on positive
  line bundles;& Proc.\ Conf.\ Complex algebraic varieties (Bayreuth,
  April~2-6, 1990), edited by K.~Hulek, T.~Peternell, M.~Schneider,
  F.~Schreyer, Lecture Notes in Math., Vol.~1507, Springer-Verlag, Berlin,
  (1992)&

\bibitem[Dem92]&J.-P.\ Demailly:& Regularization of closed positive currents 
  and Intersection Theory;& J.\ Alg.\ Geom.\ {\bf 1} (1992), 361--409&
  
\bibitem[Dem93a]&J.-P.\ Demailly:& Monge-Amp\`ere operators, Lelong
  numbers and intersection theory;& Complex Analysis and Geometry,
  Univ.\ Series in Math., edited by V.~Ancona and A.~Silva, Plenum
  Press, New-York (1993)&

\bibitem[Dem93b]&J.-P.\ Demailly:& A numerical criterion for very ample
  line bundles;& J.~Differential Geom.\ {\bf 37} (1993) 323--374&

\bibitem[Dem94]&J.-P.\ Demailly:& $L^2$ vanishing theorems for positive
  line bundles and adjunction theory;& Lecture Notes of the CIME Session
  Transcendental methods in Algebraic Geometry, Cetraro, Italy,
  July 1994, 96$\,$p, Duke e-prints {\tt alg-geom/9410022}&

\bibitem[Dol82]&I.\ Dolgachev:& Weighted projective varieties;&
  Group actions and vector fields, Proc.\ Polish-North Am.\ Semin., 
  Vancouver 1981, Springer-Verlag, Lect.\ Notes in Math.\ {\bf 956} 
  (1982) 34--71&

\bibitem[Fle89]&A.R.\ Fletcher:& Working with weighted complete
  intersections ;& Preprint MPI/89-35, Max-Planck Institut f\"ur
  Mathematik, Bonn, 1989&
  
\bibitem[FKL93]&A.\ Fujiki, R.\ Kobayashi, S.S.Y.\ Lu:& On the
  fundamental group of certain open normal surfaces;& Saitama Math.\ 
  J.\ {\bf 11} (1993), 15--20&

\bibitem[Fut83]&A.\ Futaki:& An obstruction to the existence of Einstein
  K\"ahler metrics;& Invent.\ Math.\ {\bf 73} (1983), 437--443&
  
\bibitem[Hir64]&H.\ Hironaka:& Resolution of singularities of an
  algebraic variety over a field of characteristic zero, I,II;& 
  Ann.\ Math.\ {\bf 79} (1964), 109--326&
  
\bibitem[H\"or66]&L.\ H\"ormander:& An introduction to Complex Analysis
  in several variables;& 1966, 3rd edition,
  North-Holland Math.\ Libr., Vol.~7, Amsterdam (1973)&

\bibitem[JK01]&J.M.\ Johnson, J.\ Koll\'ar:& K\"ahler-Einstein metrics on 
  log del Pezzo surfaces in weighted projective $3$-spaces;& Ann.\ Inst.\ 
  Fourier {\bf 51} (2001), 69--79&
 
\bibitem[KMM87]&Kawamata, Y., Matsuda, K., Matsuki, K:& Introduction to
  the minimal model problem;& Adv. Stud. Pure Math.\ {\bf 10}
  (1987), 283--360&
  
\bibitem [K\ampersand al92]&J.\ Koll\'ar {\rm (with 14 coauthors)}:&
  Flips and Abundance for Algebraic Threefolds;& Ast\'erisque Vol.~211
  (1992)&
  
\bibitem[Kol95a]&J.\ Koll\'ar:& Shafarevich Maps and Automorphic
  Forms;& Princeton Univ.\ Press (1995)&
  
\bibitem[Kol97]&J.\ Koll\'ar:& Singularities of pairs,
  Algebraic Geometry, Santa Cruz,  1995;& Proceedings of Symposia 
  in Pure \ Math.\ vol. 62, AMS, 1997, pages 221-287&
  
\bibitem[Lel57]&P.\ Lelong:& Int\'egration sur un ensemble analytique
  complexe;& Bull.\ Soc.\ Math.\ France {\bf 85} (1957), 239--262&
  
\bibitem[Lel69]&P.\ Lelong:& Plurisubharmonic functions and positive
  differential forms;& Gordon and Breach, New York, and Dunod, Paris
  (1969)&

\bibitem[Lic57]&A.\ Lichnerowicz:& Sur les transformations analytiques
  des vari\'et\'es k\"ahl\'eriennes;& C.~R.\ Acad.\ Sci.\ Paris
  {\bf 244} (1957), 3011-3014&

\bibitem[Lin87]&B.\ Lichtin:& An upper semicontinuity theorem for some 
leading poles of $\vert f\vert^{2s}$;& Complex analytic singularities, 
Adv.\ Stud.\ Pure Math.\ {\bf 8}, North-Holland, Amsterdam, 1987, 241--272&
  
\bibitem[Lin89]&B.\ Lichtin:& Poles of $|f(z,w)|^{2s}$ and roots of the
  $B$-function;& Ark f\"or Math.\ {\bf 27} (1989), 283--304&
  
\bibitem[Man93]&L.\ Manivel:& Un th\'eor\`eme de prolongement $L^2$ de
  sections holomorphes d'un fibr\'e vectoriel;& Math.\ Zeitschrift,
  {\bf 212} (1993) 107--122&

\bibitem[Mat57]&Y.\ Matsushima:& Sur la structure du groupe
  d'hom\'eomorphismes analytiques d'une certaine vari\'et\'e
  k\"ahl\'erienne;& Nagoya Math.\ Journal {\bf 11} (1957), 145--150&

\bibitem[Nad89]&A.M.\ Nadel:& Multiplier ideal sheaves and existence of
  K\"ahler-Einstein metrics of positive scalar curvature;& Proc.\ Nat.\ 
  Acad.\ Sci.\ U.S.A.\ {\bf 86} (1989), 7299--7300&
  
\bibitem[Nad90]&A.M.\ Nadel:& Multiplier ideal sheaves and
  K\"ahler-Einstein metrics of positive scalar curvature;& Annals of
  Math.\ {\bf 132} (1990), 549--596&
  
\bibitem[OhT87]&T.\ Ohsawa {\rm and} K.\ Takegoshi:& On the extension
  of $L^2$ holomorphic functions;& Math.\ Zeitschrift {\bf 195} (1987)
  197--204&
  
\bibitem[Ohs88]&T.\ Ohsawa:& On the extension of $L^2$ holomorphic
  functions, II;& Publ.\ RIMS, Kyoto Univ.\ {\bf 24} (1988), 265--275&

\bibitem[PS99]&D.H.\ Phong {\rm and} J.\ Sturm:& Algebraic estimates,
  stability of local zeta functions, and uniform estimates
  for distribution functions;& January 1999 preprint,
  to appear in Ann.\ of Math&

\bibitem[PS00]&D.H.\ Phong {\rm and} J.\ Sturm:& On a conjecture of Demailly 
  and Koll\'ar;& April 2000 preprint&

\bibitem[Sho92]&V.\ Shokurov:& 3-fold log flips;&
Izv.\ Russ.\ Acad.\ Nauk Ser.\ Mat.\ Vol.~56 (1992) 105--203&
  
\bibitem[Siu74]&Y.T.\ Siu:& Analyticity of sets associated to Lelong
  numbers and the extension of closed positive currents;& Invent.\ 
  Math.\ {\bf 27} (1974), 53--156&

\bibitem[Siu87]&Y.T.\ Siu:& Lectures on Hermitian-Einstein metrics for 
  stable bundles and K\"ahler-Einstein metrics;& DMV Seminar (Band 8), 
  Birkh\"auser-Verlag, Basel-Boston (1987)&

\bibitem[Siu88]&Y.T.\ Siu:& The existence of K\"ahler-Einstein metrics on
  manifolds with positive anticanonical line bundle and a suitable finite
  symmetry group;& Ann.\ of Math.\ {\bf 127} (1988), 585--627&
  
\bibitem[Siu93]&Y.T.\ Siu:& An effective Matsusaka big theorem;&
  Ann.\ Inst.\ Fourier. {\bf 43} (1993), 1387--1405&

\bibitem[Sko72]&H.\ Skoda:&Sous-ensembles analytiques d'ordre fini ou
  infini dans $\bC^n$;& Bull.\ Soc.\ Math.\ France {\bf 100} (1972),
  353--408&
  
\bibitem[Sko75]&H.\ Skoda:& Estimations $L^2$ pour l'op\'erateur
  $\ol\partial$ et applications arithm\'etiques;& S\'eminaire
  P.~Lelong (Analyse), ann\'ee 1975/76, Lecture Notes in Math.,
  Vol.~538, Springer-Verlag, Berlin (1977), 314--323&
  
\bibitem[Tia87]&G.\ Tian:& On K\"ahler-Einstein metrics on certain
  K\"ahler manifolds with $c_1(M)>0$;& Invent.\ Math.\ {\bf 89}
  (1987), 225--246&
  
\bibitem[Var82]&A.N.\ Varchenko;& Complex exponents of a singularity
  do not change along the stratum $\mu={}$constant;& Functional Anal.\ 
  Appl.\ {\bf 16} (1982), 1--9&

\bibitem[Var83]&A.N.\ Varchenko;& Semi-continuity of the complex singularity
  index;& Functional Anal.\ Appl.\ {\bf 17} (1983), 307--308&

\bibitem[Var92]&A.N.\ Varchenko;& Asymptotic Hodge structure ...;& Math.\
  USSR Izv.\ {\bf 18} (1992), 469--512&
  
\bibitem[Yau78]&S.T.\ Yau:& On the Ricci curvature of a complex
  K\"ahler manifold and the complex Monge-Amp\`ere equation~I;& Comm.\ 
  Pure and Appl.\ Math.\ {\bf 31} (1978), 339--411&

}
\vskip30pt
\noindent
(October 12, 1999; minor revision April 29, 2000; final proofs corrected on
August 24, 2001; slightly augmented on \today\ to answer a question by
JingZhou Sun)

\end